  \documentclass[12pt]{amsart}
  \usepackage{graphicx}
  \usepackage{latexsym} 
  \usepackage[all]{xy}  
  \usepackage{amsfonts} 
  \usepackage{amsthm}   
  \usepackage{amsmath}  
  \usepackage{amssymb}  
  \usepackage{enumerate}
  \usepackage[2emode]{psfrag}
  \usepackage{hyperref}
  \usepackage{color}
  \newcommand{\cC}{{C}}  

\DeclareMathOperator{\wdi}{\raisebox{0.2ex}{\scalebox{0.75}{\ensuremath{\lozenge}}}}
\DeclareMathOperator{\bdi}{\raisebox{0.2ex}{\scalebox{0.75}{\ensuremath{\blacklozenge}}}}

\DeclareMathOperator{\wt}{\raisebox{0ex}{\scalebox{0.75}{\ensuremath{\triangledown}}}}
\DeclareMathOperator{\bt}{\raisebox{0ex}{\scalebox{0.75}{\ensuremath{\blacktriangledown}}}}
\DeclareMathOperator{\wtd}{\raisebox{0.1ex}{\scalebox{0.5}{\ensuremath{\triangledown}}}}
\DeclareMathOperator{\btd}{\raisebox{0.1ex}{\scalebox{0.5}{\ensuremath{\blacktriangledown}}}}

\DeclareMathOperator{\wtb}{\overline{\raisebox{0ex}{\scalebox{0.75}{\ensuremath{\triangledown}}}}}
\DeclareMathOperator{\btb}{\overline{\raisebox{0ex}{\scalebox{0.75}{\ensuremath{\blacktriangledown}}}}}

\DeclareMathOperator{\wtu}{\raisebox{0ex}{\scalebox{0.75}{\ensuremath{\vartriangle}}}}

\def\cC{\mathcal C}
\def\cD{\mathcal D}
 
\def\can{{\rm can}} 
\def\beq{\begin{equation}} 
\def\eeq{\end{equation}}

\def\ot{{\otimes}} 
\def\W{{\bar{A}}}
\def\B{{\bar{A'}}}

\def\stac#1{\raise-.2cm\hbox{$\stackrel{\displaystyle\,\otimes\,}{\scriptscriptstyle{#1}}$}}
\def\cstac#1{\raise-.2cm\hbox{$\stackrel{\displaystyle\,\Box\,}{\scriptscriptstyle{#1}}$}}
\def\sstac#1{\otimes_{#1}}

\newtheorem{proposition}{Proposition}[section] 
  \newtheorem{lemma}[proposition]{Lemma} 
   
  \newtheorem{corollary}[proposition]{Corollary} 
  \newtheorem{theorem}[proposition]{Theorem} 

\theoremstyle{definition} 
  \newtheorem{definition}[proposition]{Definition} 
  \newtheorem{example}[proposition]{Example}

  \theoremstyle{remark} 
  \newtheorem{remark}[proposition]{Remark} 
  
\newcommand{\Section}{\setcounter{equation}{0}\section}

\newcommand{\thlabel}[1]{\label{th:#1}}
\newcommand{\thref}[1]{Theorem~\ref{th:#1}}
\newcommand{\selabel}[1]{\label{se:#1}}
\newcommand{\seref}[1]{Section~\ref{se:#1}}
\newcommand{\lelabel}[1]{\label{le:#1}}
\newcommand{\leref}[1]{Lemma~\ref{le:#1}}
\newcommand{\prlabel}[1]{\label{pr:#1}}
\newcommand{\prref}[1]{Proposition~\ref{pr:#1}}
\newcommand{\colabel}[1]{\label{co:#1}}
\newcommand{\coref}[1]{Corollary~\ref{co:#1}}
\newcommand{\relabel}[1]{\label{re:#1}}
\newcommand{\reref}[1]{Remark~\ref{re:#1}}
\newcommand{\exlabel}[1]{\label{ex:#1}}
\newcommand{\exref}[1]{Example~\ref{ex:#1}}

\newcommand{\eqlabel}[1]{\label{eq:#1}}
\newcommand{\equref}[1]{(\ref{eq:#1})}

\newcommand{\Desc}{{\rm Desc}}
\newcommand{\Hom}{{\rm Hom}}
\newcommand{\End}{{\rm End}}

\newcommand{\Ker}{{\rm Ker}\,}

\newcommand{\im}{{\rm Im}\,}

\newcommand{\Rat}{{\rm Rat}}

\def\ol{\overline}
\def\ul{\underline}

\def\ot{\otimes}

\def\CC{{\mathbb C}}

\def\MM{{\mathbb M}}
\def\NN{{\mathbb N}}

\newcommand{\Aa}{\mathcal{A}}
\newcommand{\Bb}{\mathcal{B}}
\newcommand{\Cc}{\mathcal{C}}
\newcommand{\Dd}{\mathcal{D}}

\newcommand{\Ff}{\mathcal{F}}

\newcommand{\Mm}{\mathcal{M}}

\newcommand{\sd}{{\sf d}}
\newcommand{\sm}{{\sf m}}

\def\*C{{}^*\hspace*{-1pt}{\Cc}}

\def\text#1{{\rm {\rm #1}}}

\begin{document} 

 \title{Morita theory for comodules over corings}

 \author{Gabriella B\"ohm}
 \address{Research Institute for Particle and Nuclear Physics, Budapest, 
 \newline\indent H-1525
 Budapest 114, P.O.B.\ 49, Hungary}
  \email{\href{mailto:G.Bohm@rmki.kfki.hu}{G.Bohm@rmki.kfki.hu}}
  \author{Joost Vercruysse}    
  \address{Vrije Universiteit Brussel VUB, Pleinlaan 2, B-1050,
  Brussel, Belgium} 
 \email{\href{mailto:jvercruy@vub.ac.be}{joost.vercruysse@vub.ac.be}}   
    \date{May, 2008} 
  \begin{abstract}
By a theorem due to Kato and Ohtake, any (not necessarily strict) Morita
context induces an equivalence between appropriate subcategories of
the module categories of the two rings in the Morita context. These are in
fact categories of firm modules for non-unital subrings.
We apply this result to various Morita contexts associated to a comodule
$\Sigma$ of an $A$-coring $\cC$.
This allows to extend (weak and strong) structure theorems in the literature,
in particular beyond the cases when any of the coring $\cC$ or the comodule
$\Sigma$ is finitely generated and projective as an $A$-module. That is, we
obtain relations between the category of $\cC$-comodules and the category of
firm modules for a firm ring $R$, which is an ideal of the endomorphism
algebra $\End^\cC(\Sigma)$.  
For a firmly projective comodule of a coseparable coring we prove a strong
structure theorem assuming only surjectivity of the canonical map.
\end{abstract}
  \maketitle

\Section*{Introduction}

There is a long tradition of using Morita theory in the study of Hopf-Galois
extensions and by generalization Galois corings and Galois comodules, see e.g.\
\cite{CohFishMont:Morita}, \cite{Doi:smmor}, \cite{Abuh:cleftentw},
\cite{CaDeGrVe:comcor}, \cite{BohmVer:cleftex}. 
One of the applications of Galois theory within the context of corings and
comodules is corresponding (generalized) descent theory.
That is, a study of the adjoint Hom and tensor functors between the category
of comodules over a coring and the category of modules over an appropriately
chosen algebra. In particular, finding of 
sufficient and necessary conditions for these functors to be full and
faithful. Results of this kind are referred to as (weak and strong) 
{\em structure theorems}.
If a coring $\Cc$ is finitely generated and projective as left $A$-module,
then its category of right comodules becomes isomorphic to the category of
modules over the dual ring $\*C$. Hence Galois theory for such a coring
describes in fact functors between two module categories, which explains the
relation with Morita theory. Although in general Morita contexts for comodules
can be constructed without any finiteness restriction on the coring $\Cc$, 
strictness of these Morita contexts implies such a finiteness condition
for $\Cc$, and usually as well for the comodule $\Sigma$,
  appearing in the Morita context (see
\cite[Lemma 2.5]{BohmVer:cleftex}). 
 
The standard result in Morita theory says that the connecting maps $\wt$ and
$\bt$ of a Morita context $(A,A',P,Q,\wt,\bt)$ are bijective (or equivalently
surjective, if the algebras $A$ and $A'$ have a unit) if and only if the
Morita context induces an equivalence of the module categories $\Mm_A$ and
$\Mm_{A'}$. It was a natural question to pose how distinct any Morita context
is from an equivalence of categories. This question has in fact two
categorically dual answers which were found by several authors. They say that
any Morita context induces an equivalence between certain quotient
categories of the original module categories (see \cite{Mul}) as well as an
equivalence between certain full subcategories of the original
module categories (see \cite{Kato}). The occurring full subcategories
consist of firm modules over the (possibly non-unital) rings $P\wt Q$ and
$Q\bt P$. In this paper we will extensively use this latter result due to Kato
and Ohtake \cite{Kato}. 

Firm rings and firm modules did appear in coring theory as a tool to construct
comatrix corings beyond the finitely generated and projective case. 
Recall that a (unital) bimodule (over unital rings) possesses a dual in the
bicategory of bimodules -- hence determines a comatrix coring -- if and only
if it is finitely generated and projective on the appropriate side. Without
this finiteness property it may have a dual only in a larger bicategory. In  
\cite{GomVer:ComatFirm} it was assumed that a bimodule has a dual in the
bicategory of firm bimodules over firm rings and a Galois theory in this
setting was developed.  

The aim of this paper is to merge ideas of the Kato-Ohtake Theorem on
equivalences between categories of firm modules induced by Morita contexts,
with the application of Morita theory within the framework of comodules for
corings. 
In relation with various questions, a number of Morita contexts has been
associated to a comodule. Their strictness was shown to imply (weak and
strong) structure theorems. Hereby we revisit some of these Morita contexts
and derive corresponding structure theorems by means of the Kato-Ohtake
Theorem. 
Since in this way strictness of the Morita context is no longer requested, we 
extend existing structure theorems in two ways. First,
the coring $\Cc$ will not have to be finitely generated and projective over
its base ring $A$, and secondly, the comodule $\Sigma$ will no longer have to
be finitely generated and projective over $A$. 
As a consequence, our structure theorems relate the category of
$\cC$-comodules to a category of firm modules for a firm ring (instead of
unital modules for a unital ring). That is, the framework developed in
\cite{GomVer:ComatFirm} is applied.  

Another application of our theory is to 
coseparable corings.  
It is known that coseparable corings provide a
class of examples
of firm rings (see \cite{BKW}). In particular, 
Galois theory for comodules over a coseparable coring can therefore be
reduced to Morita theory between firm rings. Applying Morita theory over firm
rings, in particular the Kato-Ohtake Theorem, to this situation, we are able
to prove 
stronger results than for arbitrary corings. Most importantly, we show that,
for a firmly projective comodule of a coseparable coring, 
surjectivity of the canonical map implies its bijectivity and this condition
is equivalent to a Strong Structure Theorem (see \thref{cosepstr}). This
theorem improves \cite[Corollary 9.4]{Ver:PhD}, \cite[5.7,
5.8]{Wis:galcom}, 
\cite[Proposition 5.6]{CaDeGrVe:comcor} and is ultimately related to
\cite[Theorem I]{Schneider90}. 
The proof of \thref{cosepstr} does not make use of any projectivity property
of the coring as a module over its base, thus it 
differs conceptually from the proofs in the papers cited above.

The paper is organized as follows. In the first section we study, and recall
facts about, general Morita theory. In \seref{firm} we collect some properties
of firm rings and study their relation with idempotent rings and
corings. \seref{Kato-Ohtake} is devoted to a full proof of the Kato-Ohtake
Theorem (which is included for the sake of completeness) and some related
results. In \seref{reduction} we develop a technique to reduce a general
Morita context to a strict Morita context over firm rings. In
\seref{firm_context} properties of Morita contexts between firm rings are
discussed. 
The results of the first section are applied to particular Morita
contexts associated to comodules in \seref{comodules}. 
The theory of \cite{GomVer:ComatFirm} 
can be applied if, for a given comodule $\Sigma$ of an $A$-coring $\Cc$, one
can find a firm ring $R$ together with a ring morphism $R\to 
\End^\Cc(\Sigma)$, such that $\Sigma$ is $R$-firmly projective as a right
$A$-module. If $\Sigma$ is a finitely generated and projective right
$A$-module, then $R$ can be 
taken equal to $\End^\Cc(\Sigma)$. In \seref{constructionR} we consider a
canonical Morita context associated to $\Sigma$ as a right $A$-module and
making use of it, we describe a situation when one can find 
such a firm ring $R$ in a general setting. 
In \cite{BohmVer:cleftex}, generalizing a construction in
\cite{CaDeGrVe:comcor}, we associated a Morita context to any $\cC$-comodule
$\Sigma$, which connects the endomorphism ring $\End^\Cc(\Sigma)$ with the
dual ring $\*C$ of $\Cc$. Its strictness was related to structure
theorems. Using the results in Section \ref{sec:Morita}, we can weaken the
assumptions made in \cite{BohmVer:cleftex}. That is, instead of assuming
surjectivity of the connecting maps, we make only assumptions on properties of
its range.  
We apply a similar philosophy to reconsider in \seref{contextextension} a
Morita context associated to a pure
coring extension in \cite{BohmVer:cleftex}.
Recall that two objects in a pre-additive category determine a Morita context
of the hom-sets. In Section \ref{sec:com_Morita} we use this method to
associate  
a canonical Morita context 
to two comodules, and show how
the structure theorems of these comodules are related. 
In particular, starting with one comodule $\Sigma$, in a favourable situation
(see Theorem \ref{thm:Morita_G_B_eq}), we associate a second comodule $B$ to
it, for which the strong structure theorem holds. Comparing the resulting
Morita context, determined by the two comodules $\Sigma$ and $B$, with the
Morita context associated to $\Sigma$ in \seref{constructionR}, we derive
structure theorems for $\Sigma$.  In the final \seref{cosep} we consider a
Strong Structure Theorem for a firmly projective comodule $\Sigma$ over a
coseparable coring $\Cc$.

{\bf Notations and conventions}
For any object $X$ in a category $\Aa$, we denote the identity morphism on $X$
again by $X$. 

Throughout the paper a {\em ring} means a module $R$ over a fixed commutative
ring $k$, together with a multiplication, i.e.\ a $k$-module map
$R \otimes_k R \to R$ satisfying the associativity constraint. 
When there is no risk of confusion, multiplication will be denoted by
juxtaposition of elements of $R$. 
In general, we do not assume that the multiplication admits a unit. In the case
when it does, i.e.\ there is an element $1_R$ in $R$ such that $r 1_R=r=1_R
r$, for all $r\in R$, then we say that $R$ is a {\em ring with unit} or a {\em
  unital ring}. A right module for a non-unital ring $R$ (over a commutative
ring $k$) is a $k$-module $M$ together with a $k$-module map $M \otimes_k R
\to M$, $m\ot r \mapsto mr$, satisfying the associativity condition
$m(rr')=(mr)r'$, for $m\in M$ and $r,r'\in R$. The category of all right
$R$-modules is denoted by $\widetilde{\Mm}_R$. By convention, for a unital
algebra $R$ we consider unital modules only. That is, right $R$-modules $M$,
such that $m 1_R=m$, for all $m\in M$. The category of unital right modules of
a unital ring $R$ is denoted by $\Mm_R$. Hom-sets in ${\widetilde \Mm}_R$ and
also in $\Mm_R$ will be denoted by $\Hom_R(-,-)$. The categories
${}_R{\widetilde{\Mm}}$ and ${}_R \Mm$ of left $R$-modules are defined
symmetrically, and hom-sets are denoted as ${}_R\Hom(-,-)$.
The categories of $R$-bimodules will be denoted by ${}_R{\widetilde{\Mm}}_R$
and ${}_R\Mm_R$, respectively, with hom-sets ${}_R\Hom_R(-,-)$. 

As in \cite{GomVer:ComatFirm}, the term {\em ideal} will be slightly abused in
the following sense. Let $\iota:R\to T$ be a morphism of (possibly non-unital) 
rings. If $R$ is a left $T$-module such that $\iota$ is left $T$-linear with
respect to this action, then we will call $R$ a left ideal for $T$, even if
$\iota$ is not necessarily injective. In particular, \cite[Lemma
  5.10]{GomVer:ComatFirm} applies to this situation. 

Let $A$ be a ring with unit. An {\em $A$-coring} is a coalgebra (comonoid) in
the monoidal category ${_A\Mm_A}$, i.e.\ a triple $(\Cc,\Delta,\varepsilon)$,
where $\Cc$ is an $A$-bimodule and the {\em coproduct} $\Delta:\Cc\to
\Cc\ot_A\Cc,\ \Delta(c)=:c^{(1)}\ot_Ac^{(2)}$ (Sweedler
notation, with implicit summation understood) and the {\em counit}
$\varepsilon:\Cc\to A$ are $A$-bimodule maps that satisfy
$(\Delta\ot_A\Cc)\circ\Delta(c)=(\Cc\ot_A\Delta)\circ\Delta(c)=:
c^{(1)} \ot_Ac^{(2)}\ot_Ac^{(3)}$ and $c^{(1)}\varepsilon(c^{(2)})= c
=\varepsilon(c^{(1)})c^{(2)}$, for all $c\in \cC$. A right 
{\em $\Cc$-comodule} consists of a right $A$-module $M$ together with a
right $A$-linear map $\rho^M:M\to M\ot_A\Cc,\ \rho^M(m)=:
m^{[0]}\ot_Am^{[1]}$ (with implicit summation), called a {\em coaction}, that
satisfies $(M\ot_A\Delta)\circ\rho^M(m)=(\rho^M\ot_A\Cc)\circ\rho^M(m)=:
m^{[0]}   
\ot_Am^{[1]}\ot_Am^{[2]}$ and $m^{[0]}\varepsilon(m^{[1]})=m$, for all $m\in
M$. For two right $\Cc$-comodules $M$ and 
$N$, a right $A$-module map $M\to N$ is said to be right {\em $\Cc$-colinear}
if $\rho^N\circ f=(f\ot\Cc)\circ \rho^M$. 
The category of all right $\Cc$-comodules and right $\Cc$-colinear maps will
be denoted as $\Mm^\Cc$ and its hom-sets will be denoted by $\Hom^\cC(-,-)$. 
The category ${}^\Cc\Mm$ of left $\Cc$-comodules, with hom-sets
${}^\cC\Hom(-,-)$, is defined symmetrically. For
the coaction $\rho^M$ on a left $\Cc$-comodule $M$ the index notation
$\rho^M(m)=m^{[-1]} \ot_A m^{[0]}$ is used, for $m\in M$.
Let $\Cc$ be an $A$-coring and $R$ any (not necessarily unital) ring. If
$M\in\Mm^\Cc$ and $M$ is a left $R$-module with multiplication map $\mu:R\ot
M\to M$ such that $\mu\in\Mm^\Cc$, then we say that $M$ is an {\em $R$-$\Cc$
bicomodule}, denoted by $M\in{_R\Mm^\Cc}$. For an extensive study of corings
and comodules we refer to the monograph \cite{BrzWis:cor}.

\Section{Morita theory} \label{sec:Morita}

\subsection{Firm modules}\selabel{firm}

In this first section of somewhat preliminary nature, we collect some facts
about firm rings and their firm modules.

Let $R$ be ring,
not necessarily having a unit. 
The Dorroh-extension of $R$ is a ring with unit: $\hat{R}=R\oplus
k$. Moreover, $\widetilde{\Mm}_R$ is isomorphic to the category $\Mm_{\hat{R}}$
of unital $\hat{R}$-modules.
The ring $R$ is a two-sided ideal in $\hat{R}$ and for all
$M\in\widetilde{\Mm}_R$ and $N\in{_R\widetilde{\Mm}}$,   
$$M\ot_RN\cong M\ot_{\hat{R}}N.$$
Let $M$ be a right $R$-module. Then the right $R$-action on $M$ induces a
right $R$-linear morphism 
$$\mu_{M,R}:M\ot_RR\to M,\quad \mu_{M,R}(m\ot_Rr)=mr.$$
Denote $MR:=\{\sum_im_ir_i~|~m_i\in M,\ r_i\in R\}$. Then obviously,
$\mu_{M,R}$ factorizes as  
\[
\xymatrix{
\mu_{M,R}:M\ot_RR \ar[r]^-{\sm_{M,R}} & MR \ar[r]^i & M ,
}
\]
where $\sm_{M,R}$ is surjective and $i$ is the obvious inclusion
map. Therefore, $MR=M$ if and only if $\mu_{M,R}$ is surjective. A ring $R$ is
said to be \emph{idempotent} if and only if $R^2:=RR=R$.  

For an arbitrary ring $R$, a right $R$-module $M$ is called
\emph{firm} if $\mu_{M,R}$ is an isomorphism. 
In this case, the inverse of $\mu_{M,R}$ will be denoted by
$$
\sd_{M,R}:M\to M\ot_RR,\quad \sd_{M,R}(m)=m^r\ot_Rr,
$$
with implicit summation understood.
The category of all firm right $R$-modules with right $R$-linear maps between
them is denoted by $\Mm_R$. (This notation is justified by the fact that a
module $M$ of a unital ring $R$ is firm if and only if it is unital.)
In the same way, we introduce the category
${_R\Mm}$ of firm left $R$-modules and left $R$-linear maps and the category
${_R\Mm_S}$ of firm bimodules where $S$ is another ring.  
Taking $M=R$, we find $\mu_R:=\mu_{M,R}=\mu_{R,M}$. Hence $R\in{\Mm_R}$ if
and only if $R\in{_R\Mm}$, i.e.\ $\mu_R$ is an isomorphism with inverse denoted
by $\sd_R$.
In this situation $R$ is called a {\em firm ring}. 
This terminology is due to Quillen \cite{Qui:firm}.  
Examples of firm rings are rings with unit, rings with local units and
coseparable corings (hence they can be constructed from
split or separable extensions of (unital) rings \cite{BKW}).  
Clearly, firm rings are idempotent, but the converse is not true. We do have,
however, the following result, extending \cite[Proposition 2.5
(1)]{Mar:Morita}. Note that, for any non-unital ring
$R$, also $R\ot_R R$ is a non-unital ring, with multiplication
\begin{equation}\label{eq:S_prod}
(r_1\ot_R r'_1)(r_2\ot_R r'_2)=r_1 r'_1 \ot_R r_2 r'_2.
\end{equation}
\begin{theorem}\thlabel{idempotentfirm}
Let $R$ be a ring (not necessarily with unit) and put $S:=R\ot_RR$ which is a
ring with multiplication  
\eqref{eq:S_prod}. If $R$ is
idempotent, then the following statements hold. 
\begin{enumerate}[(i)] 
\item If $M\in\widetilde{\Mm}_R$ such that $MR=M$, 
then $M\ot_RR\in\Mm_R$ (cf. \cite[Proposition 2.5 (1)]{Mar:Morita});
\item For $M\in\widetilde{\Mm}_R$ and $N\in{_R\widetilde{\Mm}}$, there is
  an isomorphism of $k$-modules $M\ot_RN\cong M\ot_SN$; 
\item $S$ is a firm ring;
\item The categories $\Mm_R$ and $\Mm_S$ are canonically isomorphic;
\item If $M\in\widetilde{\Mm}_R$ 
then $MR\ot_RR\in\Mm_S$;
\item For any $M\in\widetilde{\Mm}_R$, $MR\ot_R R \cong M\ot_S S$, as firm
  right $S$-modules.
\end{enumerate}
\end{theorem}

\begin{proof}
\ul{(i)}. First remark that associativity implies
$\mu_{M,R}\ot_RR=M\ot_R\mu_R$.  
Consider the following exact row in $\widetilde{\Mm}_R$.
\[
\xymatrix{
0\ar[rr] &&\Ker\mu_{M,R} \ar[rr]^i && M\ot_RR \ar[rr]^{\mu_{M,R}} &&
M \ar[rr] && 0  .
}
\]
Since the
functor $-\ot_RR:\widetilde{\Mm}_R\to\widetilde{\Mm}_R$ is right 
exact, we find the following exact row in $\widetilde{\Mm}_R$ 
\[
\xymatrix{
(\Ker\mu_{M,R})\ot_RR \ar[rr]^-{i\ot_RR} && M\ot_RR\ot_RR
\ar[rr]^-{\mu_{M,R}\ot_RR} && M\ot_RR \ar[rr] && 0 
}.
\]
If we can show that $(\Ker\mu_{M,R})\ot_RR=0$, then $\mu_{M,R}\ot_RR$ is an 
isomorphism and therefore $M\ot_RR$ is a firm right $R$-module. Take
$\sum_j m_j\ot_R r_j \ot_Rs\in (\Ker\mu_{M,R})\ot_RR$.
Since $R$ is idempotent, we can write $s=\sum_i s_is'_i\in R^2$. Hence  
\[
\sum_j m_j \ot_R r_j \ot_R s 
=\sum_{i,j} m_j\ot_R r_j \ot_R s_is'_i
=\sum_{i,j} m_j\ot_R r_j s_i\ot_R s'_i
=\sum_{i,j} m_j r_j \ot_R s_i\ot_R s'_i=\!0.
\]
Thus $(\Ker\mu_{M,R})\ot_RR=0$ as needed.\\
\ul{(ii)}. If $M\in\widetilde{\Mm}_R$, then
$M\in\widetilde{\Mm}_S$ with action
$m\cdot(r\ot_Rr')=\mu_{M,R}(m\ot_Rrr')=mrr'$, and similarly for
$N\in{_R\widetilde{\Mm}}$. 
Take $m\in M$, $n\in N$ and $r\ot_R r'\in S$, then 
$$m\cdot(r\ot_Rr')\ot_Rn=mrr'\ot_Rn=m\ot_Rrr'n=m\ot_R(r\ot_
Rr')\cdot n.$$ 
Since $R$ is idempotent, we can write any $r\in R$ as
$r=\sum_ir_ir'_i$. Therefore also
\begin{eqnarray*}
m r\ot_S n&=&\sum_i mr_ir'_i\ot_S n 
=m\cdot (\sum_ir_i\ot_Rr'_i)\ot_S n \\
&=& m \ot_S (\sum_ir_i\ot_Rr'_i)\cdot n 
= \sum_i m\ot_S r_ir'_in = m\ot_S rn .
\end{eqnarray*}
These computations show that there exist unique morphisms $j_1$ and $j_2$
which render the following diagram commutative.
\[
\xymatrix{
M\ot R\ot N \ar@<.5ex>[rr] \ar@<-.5ex>[rr] && M\ot N \ar@{=}[d] \ar[rr]^{p_1}
&& M\ot_RN \ar@<.5ex>[d]^{j_1}\\ 
M\ot S\ot N \ar@<.5ex>[rr] \ar@<-.5ex>[rr] && M\ot N  \ar[rr]^{p_2} && M\ot_
SN \ar@<.5ex>[u]^{j_2} 
}
\]
We find that $j_2\circ j_1\circ p_1=j_2\circ p_2=p_1$.
Since $p_1$ is an epimorphism, we obtain that $j_2\circ j_1=M\ot_RN$. In the
same way, $j_1\circ j_2=M\ot_SN$.\\ 
\ul{(iii)}. It follows from part (i) that $R\ot_RR\ot_RR\cong R\ot_RR$. 
Therefore $R\ot_RR\ot_RR\ot_RR\cong R\ot_RR$ as well. 
Moreover, by part (ii), $(R\ot_RR)\ot_{R\ot_RR}(R\ot_RR)\cong
(R\ot_RR)\ot_R(R\ot_RR)$. We conclude that $S\ot_SS\cong S$, i.e.\ $S$ is a
firm ring.\\ 
\ul{(iv)}. Take $M\in\Mm_R$, then $M\cong M\ot_R R$ and hence $M\ot_RR\cong
M\ot_RR\ot_RR=M\ot_RS$. By part (ii) also $M\ot_RS\cong
M\ot_SS$. Combining these isomorphisms, we find that $M\cong M\ot_SS$,
i.e.\ $M\in\Mm_S$. Conversely, if $M\in\Mm_S$, then we can define a right
$R$-action on $M$ by $m\cdot r= m^s\cdot(s r)$, where $\sd_{M,S}(m)=m^s\ot_S
s\in M\ot_SS$ is the unique element such that $m^s\cdot s= m$. Then $M$ is
firm as a right $R$-module by the following sequence of isomorphisms.
$$M\cong M\ot_SS\cong M\ot_S(S\ot_RR)\cong (M\ot_SS)\ot_RR\cong
M\ot_RR.$$ 
\ul{(v)}. This follows immediately by (i) and (iv).\\
\ul{(vi)}.  By part (i), $MR \ot_R R$ is a firm right $R$-module. 
By part (ii), $MR \ot_R S\cong MR\ot_S S$.
Since $R$ is idempotent by assumption, $MR=MS$.
$S$ is a firm ring by part (iii), hence the obvious map $MS\ot_S S \to M\ot_S
S$ has an inverse $m\ot_S ss'\mapsto ms\ot_S s'$. 
Thus the following sequence of right $S$-module isomorphisms holds. 
\[
MR \ot_R R\cong 
MR \ot_R R\ot_R R\cong
MR\ot_R S \cong
MR\ot_S S \cong
MS \ot_S S \cong M\ot_S S.
\]
\end{proof}

The following proposition provides a tool to construct idempotent rings, and
therefore firm rings in combination with the previous theorem. 

\begin{proposition}\prlabel{convolutionidempotent}
Let $\Cc$ be an $A$-coring and $T$ be an $A$-ring. If $f$ is an idempotent
element in the convolution algebra ${_A\Hom_A}(\Cc,T)$, then $\im f$ is an
idempotent ring. 
\end{proposition}

\begin{proof}
Recall that multiplication in the convolution algebra ${_A\Hom_A}(\Cc,T)$
is given by 
$$
(f*g)(c)=f(c^{(1)})g(c^{(2)}),
$$
for all $f,g\in{_A\Hom_A}(\Cc,T)$ and $c\in\Cc$. Hence $f$ is idempotent in
${_A\Hom_A}(\Cc,T)$ if and only if $f(c)=f(c^{(1)})f(c^{(2)})$, from which we
immediately deduce that $\im f$ is an idempotent ring. 
\end{proof}

\begin{example}\exlabel{convolutionidempotent}
\begin{enumerate}[(i)]
\item Let $\iota:R\to T$ be a ring morphism, where $R$ is a firm ring. We can
  regard $R$ as an $\hat{R}$-coring (see \cite[Lemma 2.1]{Ver:equi}), and
  $\iota$ makes $T$ into an $R$-ring. Multiplicativity of $\iota$
  corresponds exactly to the fact that $\iota$ is an idempotent element of the
  convolution algebra ${_{\hat{R}}\Hom_{\hat{R}}}(R,T)$. Therefore $\im\iota$
  is an idempotent ring, which can also easily be verified directly. 
\item Let $\Cc$ be an $A$-coring, then the counit
  $\varepsilon\in{_A\Hom_A}(\Cc,A)$ is 
  clearly idempotent. Hence $R=\im\varepsilon$ is an idempotent ring. In this
  situation there holds moreover a similar statement for the $\Cc$-comodules:
  for all $M\in\Mm^\Cc$, we have $M\in\Mm_R$. 
\end{enumerate}
\end{example}

\begin{remark}\relabel{functorJ}
Let $R$ be a right ideal in a unital ring $A$. Regarding $R$ as an $R$-$A$
bimodule, there is a functor 
\begin{equation}\eqlabel{functorJ}
J_R:= -\ot_R R :\Mm_R\to\Mm_A.
\end{equation}
Explicitly, for a firm right $R$-module $M$ and $m\in M$, the action by $a\in
A$ on $J_R(M)\cong M$ comes out as  
\begin{equation}\eqlabel{Aaction}
m\cdot a=m^r(ra),
\end{equation}
cf. \cite[Lemma 5.11]{GomVer:ComatFirm}.
It is straightforward to check that restricting the $A$-action on $J_R(M)$ to
$R$, we recover the original $R$-module $M$ (in particular, $J_R(M)$ is
firm as right $R$-module). 
That is to say, composing the functor $J_R:\Mm_R\to\Mm_A$ with the forgetful
functor $\Mm_A \to \widetilde{\Mm}_R$, we obtain the fully faithful inclusion
functor $\Mm_R\to \widetilde{\Mm}_R$. Thus we conclude that both $J_R$ and the
forgetful functor $\Mm_A \to \widetilde{\Mm}_R$ are fully faithful. 
\end{remark}

\begin{lemma}\lelabel{tensormodules}
Let $R$ be 
a right ideal in a (possibly non-unital) ring $A$. Then for any
$M\in{\widetilde{\Mm}}_A$ such that $MR=M$, there is a canonical isomorphism 
$$
M\ot_{R}P\cong M\ot_{A}P, \qquad \textrm{ for all
}P\in{_A{\widetilde{\Mm}}}. 
$$ 
In particular, for any $M\in{\widetilde{\Mm}}_A$, the isomorphism $M\cong
M\ot_A R$ holds if and only if $M\cong M\ot_R R$ holds.
\end{lemma}

\begin{proof}
Take $M\in{\widetilde{\Mm}}_A$ such that $MR=M$.
Since the map $\mu_{M,R}:M\ot_RR\to M$ is surjective, we find for any $m\in M$
a (not necessarily unique) element $\sum_im_i\ot_Rr_i\in M\ot_RR$ such that
$\sum_im_ir_i=m$. Therefore, 
for all $p\in P$ and $a\in A$, 
\begin{eqnarray*}
m a\ot_{R} p= \sum_im_ir_ia\ot_{R} p= \sum_im_i\ot_{R} r_iap =
\sum_im_ir_i\ot_{R} ap =m\ot_{R} ap. 
\end{eqnarray*}
Hence there exists a map $M\ot_A P \to M\ot_R P$, $m\ot_A p \mapsto m\ot_R p$,
which is easily seen to be the inverse of the epimorphism $M\ot_R P \to M\ot_A
P$, induced by the inclusion $R\to A$.  

Properties $M\cong M\ot_A R$ and $M\cong M\ot_R R$ of $M\in
{\widetilde{\Mm}}_A$ are equivalent since any of them implies that $\mu_{M,R}$ 
is surjective.  
\end{proof}

The following observation generalizes \cite[Lemma 2.1 and Theorem
  2.2]{Ver:equi}.  

\begin{theorem}\thlabel{coringideal}
An ideal $R$ in a unital ring $A$ is a firm ring if and
only if $R$ is an $A$-coring whose counit is the inclusion map $R\to A$. 
Moreover, if these equivalent conditions hold, then the category of firm right
$R$-modules is isomorphic to the category of comodules over
the $A$-coring $R$.  
\end{theorem}

\begin{proof}
Suppose first that $R$ is a firm ring
and define a coproduct $\sd_R :R\to
R\ot_RR \cong R\ot_AR$. It has a counit given by the inclusion $R\to A$.
Conversely, if $R$ is an $A$-coring with counit given by the inclusion $R\to
A$, then   
its coproduct 
$\Delta_R:R\to R\ot_A R,\ \Delta(r)=r_{(1)}\ot_Ar_{(2)}$
satisfies
$r=r_{(1)}r_{(2)}$. This implies that $\mu_R$ is surjective, i.e.\ $R$ is an
idempotent ring. Applying \leref{tensormodules} we find that
$R\ot_AR\cong R\ot_RR$ and we can easily check that $\Delta_R : R \to 
R\ot_AR\cong R\ot_RR$ is a two-sided inverse for $\mu_R$. 

Take any $M\in \Mm_A$. Using \leref{tensormodules}, under the
above conditions  
we see that
a map 
$M\to M\ot_RR\cong M\ot_AR$ is a counital coaction for the $A$-coring $R$ if
and only if it is inverse of the $R$-action $M\ot_R R \to M$. 
Therefore an $A$-module map $M \to M'$ is a morphism of firm
   $R$-modules if and only if it is a morphism of comodules.
\end{proof}

It follows by \thref{coringideal} that if $R$ is a firm ring and an ideal in
a unital ring $A$, then the functor \equref{functorJ} can be interpreted as
the forgetful functor from the category of comodules for the $A$-coring $R$ to
$\Mm_A$. Hence we obtain

\begin{corollary}\colabel{coringadjoint}
Let $R$ be a firm ring that is an ideal in a unital ring $A$. Then the
functor $J_R:\Mm_R\to\Mm_A$ 
has a right adjoint given by $-\ot_RR\simeq-\ot_AR:\Mm_A\to\Mm_R$. 
Unit and counit are given,
for all $M\in\Mm_R$ and $N\in\Mm_A$, by
\[
\sd_{M,R}: M \to J_R(M)\ot_{R}R\qquad \qquad 
\mu_{M,R}: J_R(N\ot_{R}R) \to N .
\]
Clearly $\sd_{M,R}$ is an isomorphism for all $M\in\Mm_R$, yielding another
proof of fullness and faithfulness of $J_R$ (cf. \reref{functorJ}).
\end{corollary}

\subsection{The Kato-Ohtake Theorem}\selabel{Kato-Ohtake}

In this section we prove some results concerning Morita theory for
general associative rings, with a focus on idempotent rings. This theory has
been developed in a number of papers, see e.g. \cite{Kato}, \cite{Mul}. 
Morita theory over firm rings has already been considered in
\cite{Cae:bluebook} (where firm rings are named unital rings),
\cite{Mar:Morita} and \cite{GraVit:Morita} (where firm rings are named regular
rings), however, some crucial points in the theory that will be of importance
in this note are not treated in these papers. 

Recall that a Morita context is a sextuple $(A,A',P,Q,\wt,\bt)$, consisting of
two rings $A$ and $A'$ (with or without unit), two bimodules
$P\in{_A{\widetilde{\Mm}}_{A'}}$ and  $Q\in{_{A'}{\widetilde{\Mm}}_A}$ and two
bilinear maps $\wt:P\ot_{A'}Q\to A$ and $\bt:Q\ot_AP\to A'$, that are
subjected to the following conditions 
\[
\xymatrix{
P\ot_{A'}Q\ot_{A}P \ar[rr]^-{P\ot_{A'}\bt} \ar[d]_{\wt\ot_{A}P} &&
P\ot_{A'}A' \ar[d]^-{\mu_{P,A'}}\\ 
A\ot_{A}P \ar[rr]_-{\mu_{A,P}} && P
}\qquad
\xymatrix{
Q\ot_{A}P\ot_{A'}Q \ar[rr]^-{Q\ot_{A}\wt} \ar[d]_{\bt\ot_{{A'}}Q} &&
Q\ot_{A}A \ar[d]^-{\mu_{Q,A}}\\ 
{A'}\ot_{A'}Q \ar[rr]_-{\mu_{A',Q}} &&Q\ .
}
\]
The interest in Morita contexts arises from the fact that
they can be used to study equivalences between categories. A first step is
the following well-known theorem that relates a Morita context to a pair of
functors between module categories, together with natural transformations
relating these functors. A nice formulation of this theorem makes use of the
notion of a {\em wide Morita context}, introduced in \cite{CasGT}. Let $\Aa$
and $\Bb$ be two Abelian categories, then $(F,G,\eta,\rho)$ is said to be a
right wide Morita context between $\Aa$ and $\Bb$ if and only if $F:\Aa\to
\Bb$ and $G:\Bb\to \Aa$ are right exact functors and $\eta:GF\to 1_\Aa$ and
$\rho:FG\to 1_\Bb$ are natural transformations satisfying the conditions 
$$
F\eta=\rho F \quad {\rm and} \quad G\rho=\eta G.
$$
The following lemma extends \cite[Proposition 5.2]{ChiDaNa} about Morita
contexts between unital rings.

\begin{lemma}\lelabel{MoritaWideMorita}
Let $A$ and ${A'}$ be firm rings. Then there is a bijective correspondence
between the following objects.  
\begin{enumerate}[(i)]
\item Morita contexts of the form $(A,{A'},P,Q,\wt,\bt)$, where $P\in
  {}_A\Mm_{A'}$ and $Q\in {}_{A'}\Mm_A$;
\item Right wide Morita contexts $(F,G,\omega,\beta)$ between $\Mm_A$ and
  $\Mm_{A'}$ such that $F$ and $G$ preserve direct sums; 
\item Right wide Morita contexts $(F',G',\omega',\beta')$ between ${_A\Mm}$
  and ${_{A'}\Mm}$ such that $F$ and $G$ preserve direct sums. 
\end{enumerate}
\end{lemma}

\begin{proof}
$\ul{(i)\Rightarrow (ii)}$. We can define functors $F$ and $G$ by
$F(M)=M\ot_A P$ and $G(N)=N\ot_{A'}Q$ for all $M\in\Mm_A$ and
  $N\in\Mm_{A'}$. The natural transformations $\omega$ and $\beta$ are given
  by 
\begin{eqnarray}
\xymatrix{\omega_M : M\ot_AP\ot_{A'}Q \ar[rr]^-{M\ot_A\wt} && M\ot_AA
  \ar[r]^-{\mu_{M,A}} & M} && {\rm and}\eqlabel{omega}\\
\xymatrix{\beta_N : N\ot_{A'}Q\ot_AP \ar[rr]^-{N\ot_{A'}\bt} && N\ot_{A'} {A'}
  \ar[r]^-{\mu_{N,A'}} & N} \eqlabel{beta} . 
\end{eqnarray}
$\ul{(ii)\Rightarrow(i)}$. By the Eilenberg-Watts Theorem (for the
Eilenberg-Watts Theorem over firm rings we refer to \cite{Ver:equi}), we can
write $F\simeq -\ot_AP$ with $P\in{_A\Mm_{A'}}$ and $G \simeq -\ot_{A'}Q$ with
$Q\in {_{A'}\Mm_A}$. Defining $\wt=\omega_A\circ(\sd_{A,P}\ot_{A'}Q)$ and
$\bt=\beta_{A'}\circ(\sd_{A',Q}\ot_A P)$, we easily 
       find that $(A,A',P,Q,\wt,\bt)$ is a Morita context.\\ 
The equivalence $\ul{(iii)\Leftrightarrow(i)}$ is proven symmetrically.
\end{proof}

The Kato-Ohtake Theorem says that, even
without assuming that in a Morita context $(A,{A'},P,Q,\wt,\bt)$
the rings $A$ and $A'$ are firm and their bimodules $P$ and $Q$ are firm,
there are (equivalence) 
functors $-\ot_\W P:\Mm_\W\to \Mm_\B$ and $-\ot_\B Q:\Mm_\B
\to \Mm_\W$, where $\W:={P\wt Q}$ and $\B:={Q\bt P}$ are two-sided ideals in
$A$ and $A'$, respectively. Our next task is to recall this result.
We first prove the following lemmata.  

\begin{lemma}\lelabel{Moritaadjoint}
Let $(F,G,\omega,\beta)$ be a right wide Morita context between the categories
$\Aa$ and $\Bb$. If $\omega_A$ is an isomorphism for all $A\in\Aa$ then $(F,G)$
is an adjoint pair and $F$ is a fully faithful functor. 
\end{lemma}

\begin{proof}
If $\omega:GF\to {\mathcal A}$ is a natural isomorphism then
  $\omega^{-1}:{\mathcal A}\to GF$ is the unit, while $\beta:FG\to {\mathcal
  B}$ is the counit for the adjunction $(F,G)$. 
Since the unit $\omega^{-1}$ of the adjunction is a natural isomorphism, the
left adjoint $F$ is fully faithful. 
\end{proof}

\begin{lemma}\lelabel{omegabeta}
Let $(A,{A'},P,Q,\wt,\bt)$ be a Morita context 
between not necessarily unital rings, such that the connecting map $\wt$ is
surjective. 
Then, for all $M\in{\widetilde{\Mm}}_A$, the morphism
$\omega_M$ in \equref{omega} is an isomorphism if and only if $M\in\Mm_A$. 
\end{lemma}

\begin{proof}
Suppose first that $M$ is a firm right $A$-module. We have to show that
$\omega_M$ is an isomorphism.
Since both $\mu_{M,A}$ and $\wt$ are surjective, ${\omega}_M$ is an
epimorphism. Let us prove 
that ${\omega}_M$ is also a monomorphism, i.e.\ $\Ker{\omega}_M=0$. To this
end, consider the following commutative diagram in ${\widetilde{\Mm}}_A$. 
\[
\xymatrix{
\Ker{\omega}_M \ar[rr] && M\ot_A P\ot_{A'} Q \ar[rr]^-{{\omega}_M}
&& M \ar[rr] && 0\\ 
\Ker{\omega}_M\ot_A A  \ar[u]_{\mu_{\Ker{\omega}_M,A}} \ar[rr] &&
M\ot_A P\ot_{A'} Q\ot_A A \ar[rr]^-{{\omega}_M\ot_A A}
\ar[u]^{\mu_{M\ot_A P\ot_{A'}Q,A}} && M\ot_A A \ar[u]^{\mu_{M,A}} \ar[rr]
&& 0 
}
\]
The upper row is exact as ${\omega}_M$ is an epimorphism and the exactness
of lower row follows from the fact that the functor $-\ot_A A$ is right
exact. Since $M$ is firm as a right $A$-module, $\mu_{M,A}$ is an
isomorphism. Furthermore, $\mu_{M\ot_A P\ot_{A'}Q,A}$ is surjective. Indeed,
since $\wt$ is surjective, we can find for any element $a\in A$, a (not
necessarily unique) element $\sum p_a\ot_{A'}q_a\in P\ot_{A'}Q$ such that
$\sum p_a\wt q_a=a$. Hence, for all $m\ot_A p\ot_{A'}q\in M\ot_A P\ot_{A'}Q$, 
\begin{eqnarray*}
\sum \mu_{M\ot_A P\ot_{A'}Q,A}(m^a\ot_A p_a\ot_{A'}q_a\ot_A p\wt q) 
&=&\sum m^a\ot_A p_a\ot_{A'}q_a (p\wt q)\\
&=&\sum m^a\ot_A p_a\ot_{A'}(q_a\bt p) q\\
&=&\sum m^a\ot_A p_a(q_a\bt p)\ot_{A'} q\\
&=&\sum m^a\ot_A(p_a\wt q_a) p\ot_{A'} q\\
&=&m^a\ot_A a p\ot_{A'} q\\
&=&m^aa\ot_A p\ot_{A'} q= m\ot_A p\ot_{A'} q .
\end{eqnarray*}
A diagram chasing argument shows that surjectivity of $\mu_{M\ot_A
  P\ot_{A'}Q,A}$ and injectivity of $\mu_{M,A}$ imply surjectivity
of $\mu_{\Ker{\omega}_M,A}$. Hence $\Ker{\omega}_M =(\Ker{\omega}_M)
A$. However, $(\Ker{\omega}_M)A$ contains only the zero element, as for any 
$\sum_j m_j \ot_A p_j \ot_{A'} q_j \in\Ker{\omega}_M$ and $a=\sum p_a\wt
q_a\in A$ we find that   
\begin{eqnarray*}
\sum_j m_j \ot_A p_j \ot_{A'}q_j a &=& 
\sum_{j,a} m_j\ot_A p_j \ot_{A'} q_j (p_a\wt q_a)
= \sum_{j,a} m_j (p_j \wt q_j)\ot_A p_a\ot_{A'} q_a \\
&=& \sum_{j,a} {\omega}_M(m_j\ot_A p_j\ot_{A'}q_j)\ot_A p_a\ot_{A'} q_a=0.
\end{eqnarray*}
Therefore ${\omega}_M$ is an isomorphism.

Conversely, suppose now that $\omega_M$ is an isomorphism. We need to show
that $M$ is a firm right $A$-module. For any $m\in M$, 
\begin{equation}\eqlabel{eq:omegabar}
m=({\omega}_M\circ {\omega}_M^{-1})(m)\\
= \mu_{M,A} \circ (M\ot_A\wt) \circ {\omega}_M^{-1} (m).
\end{equation}
Hence $\mu_{M,A}$ is surjective, i.e.\ $MA=M$. Then also
$(M\ot_AA)A=M\ot_AA$. Furthermore, by \equref{eq:omegabar} $\mu_{M,A}$
is a split epimorphism, proving that
$\Ker\mu_{M,A}$ is a direct summand of the right $A$-module
$M\ot_AA$. Therefore, $(\Ker\mu_{M,A})A=\Ker\mu_{M,A}$. However, for all
$\sum_j m_j\ot_A a_j\in \Ker \mu_{M,A}$ and $a'\in A$, we find that $\sum_j m_j
\ot_A a_j a'= \sum_j m_j a_j \ot_A a'=0$. So we deduce that
$\Ker\mu_{M,A}=(\Ker\mu_{M,A})A=0$. Thus $\mu_{M,A}$ is injective as well.  
\end{proof}

Symmetrically to \leref{omegabeta} one can consider a Morita context
$(A,{A'},P,Q,\wt,\bt)$ of non-unital rings, such that the connecting map $\bt$
is surjective. Then, for $N\in{\widetilde{\Mm}}_{A'}$, the morphism $\beta_N$
in \equref{beta} is an isomorphism if and only if $N\in\Mm_{A'}$.

\begin{remark}\label{rem:omega_iso}
Take a Morita context $(A,{A'},P,Q,\wt,\bt)$ of {\em unital} rings. 
\leref{omegabeta} can be applied in particular to the restricted Morita
context $(\W:=P\wt Q,A',\ol{P},\ol{Q},\wtb,\btb)$, where $\W$ is a two-sided
ideal in $A$, $\ol{P}$ is an $\W$-$A'$ bimodule and $\ol{Q}$ is an $A'$-$\W$
bimodule via the restricted $\W$-actions, $\wtb: \ol{P}\ot_{A'} \ol{Q}\to \W$
is given by corestriction of $\wt$ and $\btb:\ol{Q} \ot_\W \ol{P} \to A'$ is
equal to the composite of the epimorphism $\ol{Q}\ot_\W \ol{P}\to Q\ot_A P$
with $\bt$.
Note that, for a firm right $\W$-module $M$, we know by \leref{tensormodules}
that $M\ot_A X\cong M\ot_\W X$ for all $X\in{_A{\widetilde \Mm}}$. Therefore
the natural 
morphisms $\omega_M$ and $\ol{\omega}_M$ in \equref{omega}, corresponding to 
the original and restricted Morita contexts, are related by the following
commutative diagram.  
\[
\xymatrix{
\omega_{M}: M\ot_A P\ot_{A'} Q \ar[d]_\cong \ar[rr]^-{M\ot_A\wt} && M\ot_
A \W \ar[d]_\cong \ar[r] & M \ar@{=}[d]\\ 
\ol{\omega}_M : M\ot_\W P\ot_{A'}Q \ar[rr]^-{M\ot_\W\wtb} && M\ot_\W\W
\ar[r]^-\cong & M  .
}
\]
Thus we conclude by \leref{omegabeta} that $\omega_M$ is an isomorphism for
all $M\in\Mm_\W$. Conversely, if $\omega_M$ is an isomorphism then
$\mu_{M,\W}$ is a (split) epimorphism. Hence the vertical arrows in the
above diagram are isomorphisms by \leref{tensormodules}. Therefore also
$\ol{\omega}_M$ is an isomorphism, so $M$ is a firm $\W$-module by
\leref{omegabeta}. 

Symmetrically, for a Morita context $(A,{A'},P,Q,\wt,\bt)$ one can consider
the other restricted Morita context $(A,\B:=Q\bt P,\ol{P},\ol{Q},\wtb,\btb)$,
where $\ol{P}$ is an $A$-$\B$ bimodule and $\ol{Q}$ is an $\B$-$A$ bimodule
via the restricted $\B$-actions, $\btb: \ol{Q}\ot_{A} \ol{P}\to \B$ is given
by corestriction of $\bt$ and $\wtb:\ol{P} \ot_\B \ol{Q} \to A$ is equal to
the composite of the epimorphism $\ol{P}\ot_\B \ol{Q}\to P\ot_{A'} P$ with
$\wt$. Then the morphism $\beta_N$ in \equref{beta} is an isomorphism if and
only if $N\in\Mm_\B$. Clearly, iteration of the two constructions (in
arbitrary order) yields a Morita context 
\begin{equation}\eqlabel{eq:bar_context}
(\W,\B,\ol{P},\ol{Q},\wtb,\btb),
\end{equation}
with surjective (but not necessarily bijective) connecting maps.
\end{remark}

\begin{theorem}\cite[Theorem 2.5]{Kato} \label{thm:Kato-Ohtake}\thlabel{Kato}
Let $(A,{A'},P,Q,\wt,\bt)$ be a Morita context of unital rings and consider
the restricted Morita context 
\equref{eq:bar_context}. 
Then there is an equivalence of categories 
\[
\xymatrix{
\Mm_\W \ar@<.5ex>[rr]^-{-\ot_\W\ol{P}} && \Mm_\B
\ar@<.5ex>[ll]^-{-\ot_\B\ol{Q}} 
}.
\]
\end{theorem}

\begin{proof}
Consider the following diagram of functors
\[
\xymatrix{
\Mm_A \ar@<.5ex>[rr]^-{-\ot_A{P}} && \Mm_{A'} \ar@<.5ex>[ll]^-{-\ot_{A'}{Q}}\\
\Mm_\W \ar[u]^{J_\W} && \Mm_\B \ar[u]_{J_\B}
}
\]
where $J_\W$ and $J_\B$ are defined as in \equref{functorJ}. Recall
from \reref{functorJ} that $J_\W(M)=M$ 
as (firm) right $\W$-modules, for any $M\in\Mm_\W$.  
Hence we can apply \leref{omegabeta} to the Morita context
\equref{eq:bar_context} to conclude that
$\ol{\omega}_{J_\W(M)}=\mu_{M,\W}\circ(M\ot_\W\wtb) : M\ot_\W P\ot_{\B}Q\to M$
is  
an isomorphism 
of right $\W$-modules.  
Symmetrically, ${\overline{\beta}}_{J_\B(M')} := \mu_{M',\B}\circ (M'\otimes_\B
\btb)$ is an isomorphism, for all $M'\in \Mm_\B$.
In order to show that $J_\W(M)\ot_AP=M\ot_\W\ol{P}$ is a firm right
$\B$-module for all $M\in\Mm_\W$, we construct the inverse for the
multiplication map 
$\mu_{M\ot_\W\ol{P},\B}:M\ot_\W P\ot_\B \B \to M\ot_\W P$ as 
\[
\xymatrix{
\sd_{M\ot_\W\ol{P},\B} :
M\ot_\W P
\ar[rr]^-{\ol{\omega}^{-1}_{J_\W(M)}\ot_\W P} && M\ot_\W P\ot_\B Q\ot_\W P
\ar[rr]^-{M\ot_\W P\ot_\B \bt} 
&& M\ot_\W P\ot_\B \B
}.
\]
Thus we conclude that the functors $-\ot_\W\ol{P}=(-\ot_AP)\circ
J_\W:\Mm_\W\to\Mm_\B$ and $-\ot_\B\ol{Q}=(-\ot_{A'}Q)\circ
J_\B:\Mm_\B\to\Mm_\W$ are well-defined. Moreover,
$(-\ot_\W\ol{P},-\ot_\B\ol{Q},\ol{\omega},\ol{\beta})$ constitute a right wide
Morita context between $\Mm_\W$ and $\Mm_\B$. Since $\ol{\omega}$ and
$\ol{\beta}$ are natural isomorphisms, it follows by \leref{Moritaadjoint} 
that the context induces an equivalence of categories. 
\end{proof}

\begin{remark}
By symmetry, any Morita context $(A,{A'},P,Q,\wt,\bt)$ with restricted form
$(\W,\B,\ol{P},\ol{Q},\wtb,\btb)$ in \equref{eq:bar_context} induces as well
an equivalence of categories  
\[
\xymatrix{
{_\W\Mm} \ar@<.5ex>[rr]^-{\ol{Q}\ot_\W-} && {_\B\Mm}
\ar@<.5ex>[ll]^-{\ol{P}\ot_\B-} 
}.
\]
\end{remark}

\subsection{Reduction of a Morita context}\selabel{reduction}
Let $(A,{A'},P,Q,\wt,\bt)$ be a Morita context. In this section we extend the
construction of an associated Morita context \equref{eq:bar_context} with {\em
  surjective} connecting maps 
to appropriate (non-unital) {\em subrings} of $P\wt Q$ and $Q\bt P$.

\begin{lemma}\label{lem:B-W_prop}
Let $(A,{A'},P,Q,\wt,\bt)$ be a Morita context and let $B\subseteq Q\bt P$ be
an idempotent left ideal.
That is, assume that $BB=B$ and $Q\bt P B\subseteq B$. 
In terms of $B$, introduce the ideal $W:=PB\wt Q= P\wt BQ$ in
$A$. Consider $P$ as a $W$-$B$ bimodule and $Q$ as a $B$-$W$ bimodule via
restriction. The (non-unital) rings $B$ and $W$ obey the following
properties.
\begin{enumerate}[(i)]
\item $Q \bt PB =B$;
\item $WPB=PB$ and $BQW=BQ$; 
\item $W$ is idempotent, that is, $WW=W$;
\item $B':=QW \bt P$ is an idempotent ideal in $Q\bt P$, satisfying $B'=BQ\bt
  P$ and $W=PB'\wt Q$.
\end{enumerate}
\end{lemma}

\begin{proof}
$\underline{(i)}$ Using the assumptions that $B$ is an idempotent ring (in the
first equality) and that it is a left ideal (in the final inclusion), we
obtain a sequence of inclusions $B =BB \subseteq (Q\bt P) B \subseteq B$.\\
$\underline{(ii)}$ By construction of $W$ and part (i), associativity of
the Morita context implies $WPB =(P\wt BQ)PB= P B(Q\bt PB) = PBB = PB$. 
Symmetrically, $BQW=BQ(P\wt BQ)=B(Q\bt PB) Q=BBQ=BQ$.\\ 
$\underline{(iii)}$ Using part (ii), one deduces $WW=WPB\wt Q =PB\wt Q =W$. \\
$\underline{(iv)}$ Interchanging in part (iii) the role of $A$ with $A'$, $P$
with $Q$ and $\bt$ with $\wt$, and replacing $B$ by $W$ and $W$ by $B'$, we  
conclude that $B'$ is an idempotent ideal. Moreover,  
\[
B'=QW\bt P = Q(PB\wt Q)\bt P =(Q\bt P) B (Q\bt P) = BQ\bt P,
\]
where the last equality follows by part (i).
Since $B$ is an idempotent left ideal in $Q\bt P$, we have $B=BB \subseteq BQ
\bt P=B'$. Hence $W=PB\wt Q\subseteq PB'\wt Q$. Conversely,
\[
PB'\wt Q= PB(Q\bt P)\wt Q=(PB\wt Q)(P\wt Q)\subseteq PB\wt Q=W,
\]
since $W=PB\wt Q$ is a (right) ideal in $P\wt Q$.
\end{proof}

\begin{definition}\label{def:B-red_context}
Let $\mathbb{M}=(A,{A'},P,Q,\wt,\bt)$ be a Morita context and let $B\subseteq
Q\bt P$ be an idempotent left ideal. Introduce the ideal $W:=PB\wt
Q$ in $A$. Consider $P$ as a $W$-$B$ bimodule and $Q$ as a $B$-$W$ bimodule
via restriction. The {\em 
$B$-reduced form} of $\mathbb{M}$ is the Morita context 
\begin{equation}\label{eq:red_context}
{\overline{\mathbb{M}}}_B:=(W,B, P\ot_B B, B\ot_B Q, \ol{\wt},
\ol{\bt}),  
\end{equation}
with connecting maps
\begin{eqnarray*}
&\ol{\wt}:P\ot_B B\ot_B B\ot_B Q \to W,\qquad
&p\ot_B b \ot_B b' \ot_B q \mapsto pb \wt b' q,\nonumber \\
&\ol{\bt}:B\ot_B Q\ot_W P\ot_B B \to B,\qquad
&b \ot_B q\ot_W p\ot_B b' \mapsto b q \bt p b'.
\end{eqnarray*}
\end{definition}

One could consider many variations of the conditions on $B$, imposed in
Definition \ref{def:B-red_context}. For example, $B$
can be an ideal with respect to a ring morphism $\iota:B\to Q\bt P$. Of course
we can replace $B$ by an idempotent right ideal, or as well consider the 
$W$-reduced form of $\MM$ where $W$ is an idempotent left ideal of $P\wt
Q$. It follows from the following lemma that these approaches lead to
equivalent descriptions. 

\begin{lemma}\lelabel{conditionreduced}
Let $(A,{A'},P,Q,\wt,\bt)$ be a Morita context. Then the following statements
are equivalent. 
\begin{enumerate}[(i)]
\item There exists an idempotent left ideal $B\subseteq Q\bt P$, that is
  $BB=B$ and $Q\bt P B\subseteq B$; 
\item There exists an idempotent two-sided ideal $B'\subset Q\bt P$, that is
  $B'B'=B'$ and $Q\bt P B'\subseteq B'$ and $B'Q\bt P\subseteq B'$;  
\item There exists a firm ring $\widetilde{B}$ together with a ring morphism
  $\iota:\widetilde{B}\to Q\bt P$ such that $\widetilde{B}$ becomes a left
  ideal in $Q\bt P$, that is $\widetilde{B}$ is a left $Q\bt P$-module and
  $\iota$ is left $Q\bt P$-linear; 
\item There exists a firm ring $\widetilde{B}$ together with a ring morphism
  $\iota:\widetilde{B}\to Q\bt P$ such that $\widetilde{B}$ is a left
  $A'$-module and $\iota$ is left $A'$-linear; 
\item There exists a firm ring $\widetilde{B}'$ together with a ring morphism
  $\iota:\widetilde{B}'\to Q\bt P$ such that $\widetilde{B}'$ becomes a two
  sided ideal in $P\bt Q$, that is $\widetilde{B}'$ is a $Q\bt P$-bimodule and
  $\iota$ is $Q\bt P$-bilinear; 
\item All statements (i)-(v), where we interchange the roles of $A$ and $A'$,
  $P$ and $Q$, $\wt$ and $\bt$.  
\end{enumerate}
\end{lemma}

\begin{proof}
$\ul{(i)\Rightarrow(ii)}$. Put $B':=BQ\bt P$, as in Lemma \ref{lem:B-W_prop}
  (iv).   
\\  
$\ul{(ii)\Rightarrow(i)}$. Trivial.\\
$\ul{(i)\Rightarrow(iii)}$. We know by \thref{idempotentfirm} that
  $\widetilde{B}:=B\ot_BB$ is a firm ring. The 
  multiplication on $B$ composed by the inclusion map $B\subseteq Q\bt P$
  defines a ring map $\iota:\widetilde{B}\to Q\bt P$, which is clearly left
  $Q\bt   P$-linear.\\ 
$\ul{(iii)\Rightarrow(i)}$. Take $B=\im\iota$. Then we know by
  \exref{convolutionidempotent} (ii) that $B$ is an idempotent ring. Since
  $\iota$ is left $Q\bt P$-linear, $B$ is a left $Q\bt P$-module.\\ 
$\ul{(iii)\Rightarrow(iv)}$. Follows from the facts that $Q\bt
  P$ is a (left) ideal in $A'$ and $\widetilde{B}$ is a firm ring. Indeed,
  define an $A'$-action on $\widetilde{B}$ as $a'\cdot {\widetilde b}:=(a'
  \iota({\widetilde b}'))\cdot  {\widetilde b}\,^{\,{\widetilde b}'}$.\\
$\ul{(iv)\Rightarrow(iii)}$. Trivial.\\
$\ul{(ii)\Leftrightarrow(v)}$. Repeating the proof of
  $(i)\Leftrightarrow(iii)$, put $\widetilde{B}'=B'\ot_{B'}B'$.\\ 
$\ul{(vi)}$. Suppose $B$ exists as in $(i)$, then we know by Lemma
  \ref{lem:B-W_prop} (iii) that $W=PB\wt Q$ is an idempotent two-sided ideal
  in $P\wt Q$. This is the symmetric statement of (ii). The converse follows
  by applying the same symmetry again. 
\end{proof}

The next theorem provides us with a criterion to identify maximal ones among
idempotent rings $B$ in Lemma \ref{lem:B-W_prop}.  

\begin{theorem}\thlabel{artenian}
If $A$ is a left Artinian ring and $I$ is an ideal in $A$, then there exists a
maximal idempotent left ideal $B\subset I$. 
\end{theorem}

\begin{proof}
Put $I_1=I$. Consider $I_2\subset I_1$ as the image of the multiplication
map $\mu_1:I\ot_A I\to I$. Inductively, we define for all $n\in\NN$, $I_n$ as
the image of the multiplication map $\mu_{n-1}:I\ot_A I_{n-1}\to
I_{n-1}$. Clearly every $I_n$ is a left ideal in $I$, hence also a left ideal
in $A$. Since $A$ is left Artinian, there exists an $N\in\NN$ such that
$I_N=I_{N+1}$. 
Putting $B=I_N$, we obtain an idempotent left ideal in $I$. For any idempotent
subring $B_0$ of $I$, $B_0=B_0^N \subseteq I_N=B$. Hence $B$ is maximal in
$I$.   
\end{proof}

\begin{remark}\relabel{rem:sec_red}
In Definition \ref{def:B-red_context} we associated to a Morita context
$\mathbb{M}=(A,{A'},P,Q,\wt,$ $\bt)$ a reduced Morita context
${\overline{\mathbb{M}}}_B$ between idempotent rings with surjective
connecting maps. Using results in \thref{idempotentfirm}, one can work
equivalently with a Morita context of firm rings and their firm
bimodules. That is,  
with the same notations as in \leref{conditionreduced}, denote
$\widetilde{B}=B\ot_BB$ and $\widetilde{W}=W\ot_WW$. 
Instead of the Morita context
\eqref{eq:red_context}, one may consider  
\begin{equation}\eqlabel{2ndreduced}
(\widetilde{W},\widetilde{B},\widetilde{W} \ot_W P\ot_B \widetilde{B},
\widetilde{B}\ot_B Q \ot_W \widetilde{W},\tilde{\wt}, \tilde{\bt}) ,
\end{equation}
with connecting maps 
\begin{eqnarray*}
\tilde{\wt}&:&W\ot_W W\ot_W P \ot_B B\ot_B B\ot_B B\ot_B B \ot_B
  Q \ot_W W\ot_W W \to W\ot_W  W,\\
&&w_1  \ot_W w_2\ot_W p \ot_B b_1\ot_B b_2 \ot_B b'_1 \ot_B b'_2
  \ot_B q' \ot_W w'_1 \ot_W w'_2 \mapsto  \\
&&\qquad\qquad \qquad\qquad\qquad\qquad \qquad\qquad\qquad\qquad
w_1 w_2 (pb_1 b_2 \wt b'_1 b'_2 q')\ot_W w'_1 w'_2
, \\ 
\tilde{\bt}&:&B\ot_B B\ot_BQ\ot_W W \ot_W W\ot_W W\ot_W W \ot_W
  P\ot_B B \ot_B B\to B\ot_BB,\\
&&b_1\ot_B b_2 \ot_B q\ot_W w_1  \ot_W w_2 \ot_W w'_1  \ot_
  W w'_2\ot_W p'\ot_B b'_1 \ot_B b'_2 \mapsto \\
&&\qquad\qquad \qquad\qquad \qquad\qquad \qquad\qquad\qquad\qquad
b_1 b_2 \ot_B (qw_1 w_2 \bt w'_1 w'_2 p') b'_1 b'_2 .
\end{eqnarray*}

A few comments relating the Morita contexts \eqref{eq:red_context} and
\equref{2ndreduced} are in order.
\label{rem:B-W_eq}
\begin{enumerate}[(i)]
\item 
The 
two reduced forms \eqref{eq:red_context} and \equref{2ndreduced} of a Morita
context exist under the same equivalent conditions of
\leref{conditionreduced}.  
\item 
We know from \thref{idempotentfirm} that $\widetilde{B}$ and $\widetilde{W}$
are firm rings, $\widetilde{W} \ot_W P\ot_B \widetilde{B}$ is a firm
$\widetilde{W}$-$\widetilde{B}$ bimodule and $\widetilde{B} \ot_B Q\ot_W
\widetilde{W}$ is a firm $\widetilde{B}$-$\widetilde{W}$ bimodule.  
Therefore, since the connecting maps $\tilde{\wt}$ and $ \tilde{\bt}$ are
surjective by construction, it follows by \leref{omegabeta} that
${\widetilde{W}} \ot_{\widetilde{W}} \tilde{\wt}\cong \tilde{\wt}$ and 
${\widetilde{B}} \ot_{\widetilde{B}} \tilde{\bt}\cong \tilde{\bt}$ are
bijective, that is to say, the Morita context \equref{2ndreduced} is {\em
  strict}. 
\item
The connecting maps in both reduced Morita contexts \eqref{eq:red_context} and
\equref{2ndreduced} 
are surjective by construction. Therefore, we obtain by \thref{idempotentfirm}
and Theorem \ref{thm:Kato-Ohtake} the following commutative diagram of
category equivalences. 
\begin{equation}\eqlabel{equivalencesreduced}
\xymatrix{
\Mm_W \ar@<.5ex>[rrrrrr]^-{-\ot_W {P\ot_B B}} \ar[d]^\cong &&&&&& \Mm_B
\ar[d]^\cong 
\ar@<.5ex>[llllll]^-{-\ot_B B \ot_B {Q}\ \simeq\  -\ot_B {Q}} \\
\Mm_{\widetilde{W}} \ar@<.5ex>[rrrrrr]^-{-\ot_{\widetilde{W}} 
{\widetilde{W} \ot_W P\ot_B \widetilde{B} \simeq -\ot_W P\ot_B B}}
&&&&&& \Mm_{\widetilde{B}}
\ar@<.5ex>[llllll]^-{-\ot_{\widetilde B} \widetilde{B} \ot_B Q\ot_W
\widetilde{W} \simeq - \ot_B Q}  \ .
}
\end{equation}
That is, both reduced forms \eqref{eq:red_context} and \equref{2ndreduced} of
a Morita 
context induce (up to isomorphism of categories) the same equivalence.  
\item Let us use the notations in Lemma \ref{lem:B-W_prop} (iv),
i.e.\ put $B'=QW\bt P$. Consider the $B$-reduced
  and the $B'$-reduced forms of $\MM$, with $W=PB\wt Q=PB'\wt Q=W'$. 
Then we have equivalences
\[
\xymatrix{ 
\Mm_{W'} \ar@<.5ex>[rrrrrr]^-{-\ot_{W'} {P\ot_{B'} {B'}}} \ar@{=}[d] &&&&&&
\Mm_{B'} 
\ar@<.5ex>[llllll]^-{-\ot_{B'} {Q}} \\
\Mm_W \ar@<.5ex>[rrrrrr]^-{-\ot_W {P\ot_B B}} \
&&&&&& \Mm_B
\ar@<.5ex>[llllll]^-{-\ot_B {Q}} 
\ .}
\]
In particular, 
\begin{equation}\label{eq:B-B'}
-\ot_{B'} Q \ot_W {P\ot_B B} \simeq -\ot_{B} B:\Mm_{B'}\to \Mm_{B}
\end{equation}
is an equivalence. 
Consider now the $W$-reduced form $(W,B',W\ot_W P,Q\ot_W W, \underline{\wt},
\underline{\bt})$ of $\MM$. 
It induces equivalence functors
\[
\xymatrix{
\Mm_W \ar@<.5ex>[rrrrrr]^{-\ot_W P\simeq -\ot_W P\ot_{B'} B'} 
&&&&&& \Mm_{B'}\ar@<.5ex>[llllll]^-{-\ot_{B'} {Q}\simeq
  -\ot_{B'} {Q}\ot_W W} \ . 
}
\]
Thus we find that the $B'$-reduced and $W$-reduced Morita contexts give rise to
the same equivalences of categories. Hence also the $B$-reduced and
$W$-reduced forms give rise to the same equivalences (upto \eqref{eq:B-B'}).  
\item
Two points should be noticed here, which will be of importance later in this
paper. First, we were able to reduce our original Morita context to a strict  
Morita context \equref{2ndreduced} that induces an equivalence between
categories of firm modules over firm rings, and the functors are induced by
firm bimodules.   
Second, it is possible to represent (at least) one of the functors by the
original (possibly non-firm) bimodule from the original Morita context. 
\end{enumerate}
\end{remark}

\subsection{Morita contexts between firm rings}\selabel{firm_context}

We finish this section by extending some classical results in Morita theory
(of unital rings) to the situation of Morita contexts over firm rings, which
applies in particular to the reduced Morita context of the form
\equref{2ndreduced}.  

\begin{theorem}\thlabel{Moritafirm}
Consider a Morita context $\MM=(A,A',P,Q,\wt,\bt)$, such that $A$ is a 
firm ring and the connecting map $\wt$ is surjective. Then the following
statements hold.  
\begin{enumerate}[(i)]
\item $P\ot_{A'} Q$ is an $A$-ring with multiplication $(p\ot_{A'} q)
  (p'\ot_{A'}q'):= p(q\bt p')\ot_{A'} q'$ and unit $u:A \to P\ot_{A'} Q$
  satisfying $\wt \circ u=A$;
\item The functor $-\ot_AP:\Mm_A\to\Mm_{A'}$ is fully faithful;
\item $A\ot_A P\cong A\ot_A \ {}_{A'}\Hom(Q,A')$ and
  $Q\ot_A A \cong \Hom_{A'}(P,A')\ot_A A$,
as $A$-$A'$ bimodules and $A'$-$A$ bimodules, respectively;
\item There is a natural isomorphism $-\ot_{A'}Q \ot_A A \simeq
  \Hom_{A'}(P,-)\ot_A A$ of functors $\Mm_{A'}\to\Mm_A$ and a natural
  isomorphism $A \ot_A P \ot_{A'} -\simeq A \ot_A\ {}_{A'}\Hom(Q,-)$ of
  functors ${}_{A'}\Mm\to _A \Mm$;
\item $A$ is a left ideal in $\End_{A'}(P)$ and a right ideal in
  ${}_{A'}\End(Q)^{op}$; 
\item $Q\ot_A A$ is a generator in $\Mm_A$ and $A\ot_A P$ is a generator in
  ${}_A \Mm$; 
\item If in addition $P$ is a firm left $A$-module or $Q$ is a firm right
  $A$-module, then $\wt$ is bijective.
\end{enumerate}
\end{theorem}

\begin{proof}
\ul{(i)}. 
Since $A$ is a firm ring, it follows by \leref{omegabeta} that $A\ot_A 
\wt:A\ot_A P \ot_{A'} Q \to A \ot_A A$ is bijective. Hence there is an
$A$-bimodule map 
\[
u:= (\mu_{A,P}\ot_{A'} Q)\circ (A\ot_A \wt)^{-1} \circ \sd_A:A \to P\ot_
{A'} Q,\qquad a\mapsto p_a\ot_{A'} q_a,
\]
with implicit summation understood. Since $\wt$ is left $A$-linear,
\[
\wt\circ u =
\wt \circ (\mu_{A,P}\ot_{A'} Q)\circ (A\ot_A \wt)^{-1} \circ \sd_A=
\mu_A \circ (A\ot_A \wt)\circ (A\ot_A \wt)^{-1} \circ \sd_A=A.
\]
Moreover, for $p\in P$, $q\in Q$ and $a\in A$,
\begin{eqnarray*}
&&p(q\bt p_a)\ot_{A'} q_a=p\ot_{A'} q (p_a\wt q_a) = p\ot_{A'} qa
\quad\textrm{and}\quad\\
&&p_a(q_a\bt p)\ot_{A'} q=(p_a\wt q_a) p\ot_{A'} q = ap\ot_{A'} q.
\end{eqnarray*}
Thus $P\ot_{A'} Q$ is an $A$-ring with the stated product and unit $u$.

\ul{(ii)}. Follows directly form \leref{Moritaadjoint} and
  \leref{omegabeta}. 

\ul{(iii)}. Consider the $A$-$A'$ bimodule map
\[
\alpha: A\ot_A P \to A \ot_A\ {}_{A'} \Hom(Q,A'),\qquad a\ot_A p \mapsto
a\ot_A (-\bt p).
\]
In terms of the map $u:a \mapsto p_a \ot_{A'} q_a$ in part (i), its inverse is
given by  
\[
\alpha^{-1}: A \ot_A\ {}_{A'} \Hom(Q,A')\to A\ot_A P,\qquad a_1 a_2 \ot_
A \varphi \mapsto a_1 \ot_A p_{a_2}\varphi(q_{a_2}).
\]
The other isomorphism $Q\ot_A A \cong \Hom_{A'}(P,A')\ot_A A$ follows by
symmetrical reasoning.

\ul{(iv)}. A natural transformation $-\ot_{A'}Q\ot_A A \to
\Hom_{A'}(P,-)\ot_A A$ is given, for $M\in \Mm_{A'}$, by the right $A$-module
map
\[
\Phi_M:M \ot_{A'}Q\ot_A A \to \Hom_{A'}(P,M)\ot_A A, \qquad 
m\ot_{A'} q\ot_A a \mapsto m(q\bt -)\ot_A a.
\]
In terms of the map $u:a \mapsto p_a \ot_{A'} q_a$ in part (i), its inverse is
given by  
\[
\Phi_M^{-1}:\Hom_{A'}(P,M)\ot_A A \to M \ot_{A'}Q\ot_A A, \qquad 
\varphi\ot_A a_1 a_2 \mapsto \varphi(p_{a_1})\ot_{A'} q_{a_1}\ot_A a_2.
\]
The other natural isomorphism $A \ot_A P \ot_{A'} -\simeq A \ot_A\
{}_{A'}\Hom(Q,-)$ is proven symmetrically.

\ul{(v)}. As in the proof of part (i), \leref{omegabeta} implies that $\wt
\ot_A A: P\ot_{A'} Q \ot_A A \to A \ot_A A\cong A$ is an
isomorphism. Moreover, by part (iv), $P\ot_{A'} Q \ot_A A \cong
\End_{A'}(P)\ot_A A$. By \cite[Lemma 5.10]{GomVer:ComatFirm}, the combined
isomorphism $\End_{A'}(P)\ot_A A \cong A$  means exactly that $A$ is a left
ideal in $\End_{A'}(P)$. The other claim follows symmetrically. 

\ul{(vi)}. To any $p\in P$ we can associate a map
$f^p\in\Hom_A(Q\ot_A A ,A)$, defined by $f^p(q\ot_A a):=(p\wt q)a$. 
In terms of the map $u:a \mapsto p_a \ot_{A'} q_a$ in part (i), 
$f^{p_a}(q_a \ot_A {\tilde a})=a{\tilde a}$, for all $a,{\tilde a}\in
A$. Hence it follows by
the firm property of $A$ that the evaluation map $\Hom_A(Q\ot_A A ,A)\ot_{A'}
Q\ot_A A\to A$ is surjective. Consider the following commutative
diagram with obvious maps, for all $M\in\Mm_A$.
\[
\xymatrix{
M\ot_A\Hom_A(Q\ot_A A,A)\ot_{A'} Q\ot_A A \ar[rr] \ar[d] && M\ot_A A
\ar[d]^\cong \\ 
\Hom_A(Q\ot_A A,M)\ot_{A'} Q\ot_A A\ar[rr] && M
}
\]
Since the map in the upper row is an epimorphism, we find that the map in the
lower row is an epimorphism as well, i.e.\ $Q\ot_A A$ is a generator for
$\Mm_A$. It follows by a symmetrical reasoning that $A\ot_A P$ is a generator
in ${}_A \Mm$.

\ul{(vii)}. Assume that $P$ is a firm left $A$-module. Since $\wt$ is
surjective, it follows from \leref{omegabeta} that $\omega_M$ is an
isomorphism for all $M\in\Mm_A$, so for $M=A$. Since $\omega_A$ and $\wt$
differ by isomorphisms ($\mu_A$ and $\mu_{A,P}\ot_{A'} Q$), $\wt$ is an
isomorphism, too. The case when $Q$ is a firm right $A$-module is treated
symmetrically. 
\end{proof}

Generalizing finitely generated projective modules over unital rings to firm
modules over firm rings, the notion of {firm projectivity} was introduced
in \cite{Ver:equi}. Since this notion plays a central role also in the present
paper, in the next theorem we recall some facts about it (without proof).

\begin{theorem}\cite[Theorem 2.4]{Ver:equi},\cite[Theorem 2.51]{Ver:PhD}
\label{thm:firm_proj}
Let $R$ and $S$ be firm rings and $\Sigma$ a firm $R$-$S$ bimodule. The
following 
statements are equivalent. 
\begin{enumerate}[(i)]
\item $\Sigma$ possesses a right dual (equal to $\Sigma^*\ot_R R$) in the
  bicategory of firm bimodules (formulated sometimes as $\Sigma$ is a
  connecting bimodule in a {\em comatrix coring context}); 
\item There is a natural isomorphism $\Hom_S(\Sigma,-)\ot_R R \simeq -\ot_S
  \Sigma^* \ot_R R$ of functors $\Mm_S \to \Mm_R$;
\item There is a non-unital ring map $R\to \Sigma\ot_S \Sigma^*$, which
  induces the original left $R$-action on $\Sigma$.
\end{enumerate}
Here $\Sigma^*:=\Hom_S(\Sigma,S)$ and a (non-unital) multiplication in
$\Sigma\ot_S \Sigma^*$ is induced by the evaluation map. A bimodule $\Sigma$
obeying these equivalent properties is said to be an {\em $R$-firmly
projective right $S$-module}. 
\end{theorem}

\begin{corollary}\label{cor:firm_ideal}
Let $\MM=(A,A',P,Q,\wt,\bt)$ be a Morita context with unital rings $A$ and
$A'$ and let $R$ be a firm ring and  left ideal in $P\wt Q$. Then the following
assertions hold. 
\begin{enumerate}[(i)]
\item There exists an $R$-bimodule map $\wtu:R\to P\ot_{A'} Q$, such that
  $(\wt\circ \wtu)(r)=r$, for all $r\in R$;
\item $R\ot_R P$ is an $R$-firmly projective right $A'$-module.
\end{enumerate}
\end{corollary}

\begin{proof}
The $R$-reduced Morita context $(R,S:=PR\wt Q,R\ot_R P,Q\ot_R R,\wtb,\btb)$
satisfies all assumptions in \thref{Moritafirm}. Hence by \thref{Moritafirm}
(i), there is an $R$-bimodule map $\ol{u}: R \to R\ot_R P\ot_S Q \ot_R R$,
such that $\wtb \circ \ol{u}=R$.

\ul{(i)}. 
In terms of $\ol{u}$, introduce the composite map
\[
\wtu:
\xymatrix{
R\ar[r]^-{\ol{u}}&
R\ot_R P\ot_S Q \ot_R R \ar[rrr]^-{\mu_{R,P}\ot_S  \mu_{Q,R}}&&&
P\ot_S Q \ar[r]&
P\ot_{A'} Q
},
\]
where the rightmost arrow denotes the canonical epimorphism. The map $\wtu$ is
$R$-bilinear and satisfies $(\wt\circ \wtu)(r) =(\wtb\circ \ol{u})(r)=r$, for
all $r\in R$.

\ul{(ii)}.
$R\ot_R P$ is a firm $R$-$A'$ bimodule hence the claim is proven by
construction of a (non-unital) ring map $\jmath':R\to (R\ot_R P)\ot_{A'}
(R\ot_R P)^*$, which is compatible with the left $R$-action on $R\ot_R
P$. Consider the left $A'$-module map   
\[
f:Q\to 
(R\ot_R P)^*, \qquad q\mapsto \big(\ r\ot_R p \mapsto q\bt
rp\ \big). 
\]
It can be used to construct an $R$-bimodule map 
\[
\jmath':= (R\ot_R P\ot_{A'} f) \circ (R\ot_R \wtu)\circ \sd_R:R \to
(R\ot_R P)\ot_{A'} 
(R\ot_R P)^*.
\]
A straightforward computation yields $\jmath'(r_1)\jmath'(r_2)=r_1
\jmath'(r_2)$ for all $r_1,r_2\in R$ hence, by its left $R$-linearity,
multiplicativity of $\jmath'$. Its compatibility with the left $R$-action
$\mu_{R}\ot_R P$ follows immediately by associativity of a Morita context and
part (i).
\end{proof}

\Section{Applications to comodules over a coring}\selabel{comodules}

In this section we apply the theory developed in Section \ref{sec:Morita} to
various Morita contexts associated to comodules of corings. 
Let $\Cc$ be a coring over a unital ring $A$ and let $\Sigma$ be a right
$\Cc$-comodule. Denote $T=\End^\Cc(\Sigma)$, then there exists a pair of
adjoint functors   
\begin{equation}\label{eq:T_adj}
\xymatrix{
\Mm_T \ar@<.5ex>[rrr]^-{-\ot_T\Sigma} &&& \Mm^\Cc
\ar@<.5ex>[lll]^-{\Hom^\Cc(\Sigma,-)} 
}
\end{equation}
whose unit and counit are given by
\begin{eqnarray}
\label{eq:T-unit}
\nu_N:N\to\Hom^\Cc(\Sigma,N\ot_T\Sigma),&&\nu_N(n)(x)=n\ot_Tx;\\
\label{eq:T-counit}
\zeta_M:\Hom^\Cc(\Sigma,M)\ot_T\Sigma\to M,&&\zeta_M(\varphi\ot_Tx)=\varphi(x);
\end{eqnarray}
for all $N\in\Mm_T$ and $M\in\Mm^\Cc$.
Galois theory for the comodule $\Sigma$ 
includes the study of
this pair of adjoint functors, in particular it concerns the question whether
these functors or their (co)restrictions are fully faithful. 

\subsection{Morita contexts associated to a comodule}\selabel{constructionR}

A first Morita context can be associated to a comodule $\Sigma$ of an
$A$-coring $\Cc$ by considering it as right $A$-module. Then we can associate
to $\Sigma$ and $\Sigma^*:=\Hom_A(\Sigma,A)$ a Morita context 
\begin{equation}\label{eq:A_mod_context}
(S=\End_A(\Sigma),A,\Sigma,\Sigma^*,\wt,\bt),
\end{equation}
 as in
\cite[Section II.4]{Bass}. The connecting maps are in this situation given by 
\begin{eqnarray}
\wt:\Sigma\ot_A\Sigma^* \to S, && x\wt\xi=x\xi(-);
\eqlabel{eq:S_context} \\  
\bt:\Sigma^*\ot_{S}\Sigma \to A, && \xi\bt x=\xi(x).\nonumber
\end{eqnarray}
If we put $\bar{S}=\Sigma\wt\Sigma^*$, then we can restrict our Morita
context to $(\bar{S},A,\Sigma,\Sigma^*,\bar{\wt},\btb)$, where we regard the
restricted actions on $\Sigma$ and $\Sigma^*$, the corestriction $\wtb$ of
$\wt$ and the composite $\btb$ of the canonical epimorphism with $\bt$. Since
$\bar{\wt}$ is surjective by construction, we obtain by \leref{Moritaadjoint}
and \leref{omegabeta} an adjunction 
\begin{equation}\eqlabel{adjunction1}
\xymatrix{
\Mm_{\bar{S}} \ar@<.5ex>[rr]^-{-\ot_{\bar{S}}{\Sigma}} && \Mm_{{A}}
\ar@<.5ex>[ll]^-{-\ot_A {\Sigma}^*} 
} .
\end{equation}

\begin{lemma}\lelabel{ringextensionadjoint}
Let $\iota:R\to S$ be a ring morphism where $S$ is any (possibly non-unital,
possibly non-firm) ring and $R$ is a firm ring. Then the functor
$-\ot_RS:\Mm_R\to\Mm_S$ has a right adjoint given by $-\ot_RR:\Mm_S\to\Mm_R$. 
\end{lemma}

\begin{proof}
Consider $S$ as an $R$-bimodule with actions given by $r\cdot s \cdot r'=
\iota(r)s\iota(r')$ for all $r,r'\in R$ and $s\in S$.  Let us first check
that the functor $-\ot_RS:\Mm_R\to\Mm_S$ is well-defined. Take any
$M\in\Mm_R$, then 
$M\ot_RS$ is a firm right $S$-module, that is, the multiplication map  
$$
\mu:M\ot_RS\ot_SS\to M\ot_RS,\qquad \mu(m\ot_Rs\ot_St)=m\ot_Rst;
$$ 
has a two-sided inverse
$$
\sd:M\ot_RS\to M\ot_RS\ot_SS,\qquad  \sd(m\ot_Rs)=m^r\ot_R\iota(r)\ot_Ss.
$$
Finally, let us give the unit $\alpha$ and counit $\beta$ for the adjunction,
and leave other verifications to the reader. 
\begin{eqnarray*}
&\alpha_M:M\to M\ot_RS\ot_RR,\qquad& \alpha_M(m)=m^r\ot_R\iota(r^t)\ot_Rt;\\
&\beta_N:N\ot_RR\ot_RS \to N, \qquad & \beta(n\ot_R r\ot_R s)= n\cdot
  \iota(r)s , 
\end{eqnarray*}
for all $M\in\Mm_R$ and $N\in\Mm_S$. 
\end{proof}

\begin{proposition}\prlabel{ConstructionR}
Let $\Cc$ be a coring over a unital ring $A$ and $\Sigma$ be a
right $\Cc$-comodule. Consider the 
Morita context \eqref{eq:A_mod_context} associated to $\Sigma$ as right
$A$-module 
and put $\bar{S}=\Sigma\wt\Sigma^*$ as before. 
If there exists a firm ring $R$ together with a ring morphism $\iota:R\to
\bar{S}$, then $R\ot_R\Sigma$ is $R$-firmly projective as a right
$A$-module. If there exists moreover a ring morphism $\iota':R\to
T=\End^\Cc(\Sigma)$ then $\Sigma$ (and therefore $R\ot_R\Sigma$ as well) is an
$R$-$\Cc$ bicomodule.  
\end{proposition}

\begin{proof}
Combining the adjunction \equref{adjunction1} with the adjunction in
\leref{ringextensionadjoint}, we find that the functor
$-\ot_R {\bar{S}}\ot_{\bar{S}}{\Sigma}:\Mm_R\to\Mm_A$ has a right adjoint
given by $-\ot_A\Sigma^*\ot_RR$, cf. 
\[
\xymatrix{
\Mm_R \ar@<.5ex>[rr]^-{-\ot_R\bar{S}} 
&& \Mm_{\bar{S}} \ar@<.5ex>[rr]^-{-\ot_{\bar{S}}{\Sigma}}
\ar@<.5ex>[ll]^-{-\ot_RR} && \Mm_{{A}}
\ar@<.5ex>[ll]^-{-\ot_{A}{\Sigma}^*} 
} .
\]
On the other hand, there is a natural isomorphism $-\ot_R
{\bar{S}}\ot_{\bar{S}}{\Sigma} \simeq -\ot_R\Sigma$,
given for $M\in \Mm_R$ by
\[
m\ot_R s \ot_{\bar{S}} x \mapsto m \ot_R sx,\qquad \textrm{\and}\qquad 
m\ot_R x \mapsto m^r \ot_R \iota(r)\ot_{\bar{S}} x.
\]
Using the characterization of $R$-firmly projective modules in
Theorem \ref{thm:firm_proj}, we find that
$R\ot_R\Sigma$ is $R$-firmly projective as a right $A$-module. 

Clearly, $\Sigma\in{_R\Mm^\Cc}$ if and only if the left action of $R$ on
$\Sigma$ induces a ring morphism $R\to \End^\Cc(\Sigma)$. 
\end{proof}

Consider a coring $\cC$ over a unital ring $A$ and a firm ring $R$. Any
$R$-$A$-bimodule $\Sigma$ determines an adjunction 
\[
\xymatrix{
\Mm_R \ar@<.5ex>[rrrrrr]^-{-\ot_R\Sigma}
&&&&&& 
\Mm_A
\ar@<.5ex>[llllll]^-{\Hom_A(\Sigma,-)\ot_RR 
}  
} .
\]
Since replacing the $R$-$A$ bimodule $\Sigma$ by $R\ot_R \Sigma$ we obtain
naturally isomorphic functors $-\ot_R\Sigma\simeq-\ot_RR\ot_R\Sigma:\Mm_R \to
\Mm_A$, hence also the right adjoints are naturally isomorphic by
\begin{equation}\eqlabel{chi}
\chi_M:\Hom_A(R\ot_R\Sigma,M)\ot_RR\to\Hom_A(\Sigma,M)\ot_RR.
\end{equation}
If in addition $\Sigma$ is an
$R$-$\Cc$ bicomodule, 
then 
we have the following pair of adjoint functors 
\begin{equation}\eqlabel{HomTensorSigma}
\xymatrix{
\Mm_R \ar@<.5ex>[rrrrrr]^-{F_\Sigma=-\ot_R\Sigma\simeq-\ot_RR\ot_R\Sigma}
&&&&&& 
\Mm^\Cc
\ar@<.5ex>[llllll]^-{G_\Sigma=\Hom^\Cc(R\ot_R\Sigma,-)\ot_RR\simeq
  \Hom^\Cc(\Sigma,-)\ot_RR}  
} .
\end{equation}
If $R$ is equal to the unital ring $\End^\cC(\Sigma)$ then
\equref{HomTensorSigma} reduces to the adjunction \eqref{eq:T_adj}. 
Since $M\ot_RR\ot_R\Sigma\cong M\ot_R\Sigma$ as right $\cC$-comodules for all
$M\in\Mm_R$, we find that 
the upper functors are naturally isomorphic, indeed. The natural isomorphism
between the lower functors follows from the uniqueness of the right adjoint.
Unit and counit of the adjunction \equref{HomTensorSigma} are given explicitly
by 
\begin{eqnarray}\eqlabel{unitcounit}
&\nu_N: N\to \Hom^\Cc(\Sigma,N\ot_R\Sigma)\ot_RR,\qquad & \nu_N(n)=(x\mapsto
n^r\ot_Rx)\ot_R r,\\ 
&\zeta_M: \Hom^\Cc(\Sigma,M)\ot_R\Sigma\to M,\qquad \qquad & \zeta_M(f\ot_R
x)= f(x),
\end{eqnarray}
for $N\in \Mm_R$ and $M\in \Mm^\cC$.

In any case when $R\ot_R \Sigma$ is an $R$-$\cC$
bicomodule that is $R$-firmly projective as a right $A$-module, the
theory developed in \cite{GomVer:ComatFirm} can be applied to it. (We refer to
\cite[Section 4.2.5]{Ver:PhD} for a more detailed treatment of structure
theorems.)  
The above observations make it possible to translate occurring properties of
$R\ot_R \Sigma$ to properties of $\Sigma$. Thus we obtain following 

\begin{corollary}\colabel{FirmGalois}
Let $\Cc$ be a coring over a unital ring $A$ and $\Sigma$ be a
right $\Cc$-comodule. Consider the Morita context \eqref{eq:A_mod_context}
associated to $\Sigma$ as right $A$-module 
and put $\bar{S}=\Sigma\wt\Sigma^*$ as before. 
Assume that there exists a firm ring $R$ together with ring morphisms
$\iota:R\to \bar{S}$
and $\iota':R\to T=\End^\Cc(\Sigma)$. 
Then the following statements hold.
\begin{enumerate}[(i)]
\item $\Sigma^\dagger:=\Sigma^*\ot_RR$ is a $\Cc$-$R$ bicomodule;
\item There exists a comatrix coring $\Sigma^\dagger\ot_R\Sigma$ over $A$;
\item The map $\can:\Sigma^\dagger\ot_R\Sigma\to\Cc,\
  \can(\xi\ot_Rr\ot_Rx)=\xi(rx^{[0]})x^{[1]}$ is an $A$-coring morphism; 
\item The inner and outer triangles of the following diagram of adjoint
  functors commute (upto natural isomorphism)
\[
\xymatrix{ 
\Mm_R \ar@<-0.5ex>[rrrr]_-{-\ot_R\Sigma} \ar@<0.5ex>[rrrrdd]^-{-\ot_R\Sigma}
&&&& \Mm_A \ar@<-0.5ex>[llll]_-{-\ot_A\Sigma^\dagger}
\ar@<0.5ex>[dd]^-{-\ot_A\Cc}\\ 
\\
&&&& \Mm^\Cc \ar@<0.5ex>[uu]^-{\Ff^\Cc}
\ar@<0.5ex>[lllluu]^-{\Hom^\Cc(\Sigma,-)\ot_RR\hspace{.5cm}} \ ,
}
\]
where $\Ff^\Cc$ denotes the forgetful functor;
\item There is an adjunction $(-\ot_R \Sigma, -\ot^\Cc
  \Sigma^\dagger)$, where $\ot^\Cc$ denotes cotensor product over $\cC$. Unit
  and counit of the adjunction are, for $M\in \Mm^\Cc$  and $N\in \Mm_R$, 
\begin{eqnarray*}
&N \to (N\ot_R \Sigma)\ot^\Cc \Sigma^\dagger, \qquad 
&n \mapsto (n^r)^{r'}\ot_R e_{r'}\ot_A (f_{r'}\ot_R r),\\
&(M\ot^\Cc \Sigma^\dagger)\ot_R \Sigma \to M, \qquad
&m\ot_A (f \ot_R r) \ot_R x \mapsto m f(rx) ,
\end{eqnarray*}
where the map $(\mu_{R,\Sigma} \ot_A \mu_{\Sigma^*,R})\circ (R\ot_R
  \Sigma\ot_A \chi_A)\circ (\jmath\ot_R R)\circ \sd_R: R \to
  \Sigma\ot_A \Sigma^*$, $r\mapsto e_r\ot_A f_r$ is obtained from the
  (non-unital) ring morphism $\jmath: R\to (R\ot_R \Sigma)\ot_A (R\ot_R
  \Sigma)^*$,  coming from 
  firm projectivity of $R\ot_R \Sigma$, and the isomorphism $\chi_A: (R\ot_R
  \Sigma)^\dagger\to \Sigma^\dagger$, $\varphi\ot_R rr' \mapsto \varphi(r\ot_R
  -)\ot_R r'$ in \equref{chi}. 

By uniqueness of a right adjoint, there is a natural isomorphism 
$-\ot^\Cc \Sigma^\dagger\simeq \Hom^\cC(\Sigma,-) \ot_R R$; 

\item If the functor $\Hom^\Cc(\Sigma,-)\ot_RR:\Mm^\Cc\to\Mm_R$ is fully
  faithful then $\can$ is an isomorphism of $A$-corings;

\item The functor $\Hom^\Cc(\Sigma,-)\ot_RR:\Mm^\Cc\to\Mm_R$ is fully faithful
  and $\Cc$ is flat as a left $A$-module if and only if $\can$ is an
  isomorphism and $R\ot_R\Sigma$ is flat as a left $R$-module (meaning that
  the functor $-\otimes_R R\ot_R\Sigma$, from the category $\Mm_{\hat{R}}$ of 
  modules of the Dorroh-extension ${\hat{R}}$ to the category $Ab$ of Abelian
  groups, is (left) exact);  

\item If $\Sigma$ is totally faithful as a left $R$-module (meaning that, for
  any $N\in \Mm_R$, $N\ot_R \Sigma=0$ implies $N=0$), then the functor
  $-\ot_R\Sigma:\Mm_R\to\Mm^\Cc$ is fully faithful. The converse holds if the
  map $\Sigma\ot_A\Sigma^\dagger\to\End_A(\Sigma)$ is a pure left
  $R$-module monomorphism; 
\item 
If $\Cc$ is flat as a left $A$-module, then $-\ot_R\Sigma:\Mm_R\to\Mm^\Cc$ is
  an equivalence if and only if $\Sigma$ is faithfully flat as left $R$-module
  and $\can$ is an isomorphism of $A$-corings. Moreover, in this situation $R$
  is a left ideal in $\End^\Cc(\Sigma)$. 
\end{enumerate}
\end{corollary}

\begin{example}
Let $\Cc$ be a coring over a unital ring $A$ and let $\Sigma$ be a right
$\cC$-comodule.  
In this example we provide an explicit construction of the firm ring $R$ in
\prref{ConstructionR} in appropriate situations. Let
$\bar{S}:=\Sigma\wt\Sigma^*$ be defined in terms of the connecting map
\equref{eq:S_context}, and $T:=\End^\Cc(\Sigma)$. Put $B:=\bar{S}\cap T$. Since
$\bar{S}$ is an ideal in $S=\End_A(\Sigma)$ by 
construction and $T\subset S$, we conclude that $B$ is an ideal in $T$. If $T$
is a (left) Artinian ring, then we can apply \thref{artenian} to obtain an
idempotent (left) ideal $B'\subset B$, which is still a (left) ideal in $T$,
and a subring of $\bar{S}$. If we put now $R:=B'\ot_{B'}B'$, then $R$ is a firm
ring and the ring morphisms $\iota: R=B'\ot_{B'}B'\to \bar{S}$ and $\iota':
R=B'\ot_{B'}B'\to T$ are given by the multiplication on $B'$. Remark that $R$
is still a left $T$-module. 
\end{example}

In a recent paper \cite{BohmVer:cleftex} we associated also another Morita
context to a comodule. Consider a coring $\cC$ over
a unital ring $A$ and a right $\cC$-comodule $\Sigma$. There exists a Morita
context 
\begin{equation}\eqlabel{contextSigma}
\MM(\Sigma)=(T,\*C,\Sigma,Q,\wt,\bt), 
\end{equation}
connecting the unital rings $T := \End^\cC(\Sigma)$ and
${}^*\cC={}_A\Hom(\cC,A)$. The bimodule $Q$ is defined by 
\begin{eqnarray}
\label{eq:Q}
Q
&=&
\{\ q\in \Hom_A(\Sigma,{}^*\cC)\ |
\ \forall x\in \Sigma,c\in \cC\quad
q(x^{[0]})(c)x^{[1]}=c^{(1)}q(x)(c^{(2)})\ \}\\
&\cong&\{\ q\in {}_A{\rm Hom}(\cC,\Sigma^*)\ |\
\forall x\in \Sigma, c\in \cC\quad 
c^{(1)}q(c^{(2)})(x)=q(c)(x^{^{[0]}})x^{[1]}\ \}.\nonumber
\end{eqnarray}
The two forms of $Q$ are related by interchanging the order of the arguments
and their parallel use should cause no confusion. The connecting maps are
\begin{eqnarray}
\bt:Q\ot_{T} \Sigma &\to {}^*\cC,\qquad q\ot_{T} x &\mapsto
q(x),\label{eq:bt}\\
\wt:\Sigma\ot_{{}^*\cC} Q&\to T,\qquad x\ot_{{}^*\cC} q&\mapsto
xq(-).\label{eq:wt}
\end{eqnarray}
For more details we refer to \cite[Section 2]{BohmVer:cleftex}.

\begin{theorem}
For a unital ring $A$,
let $\Cc$ be an $A$-coring and $\Sigma$ a right $\Cc$-comodule. Consider the
Morita context \equref{contextSigma} associated to $\Sigma$. If there exists a
firm ring $R$ together with a ring morphism $\iota:R\to \Sigma\wt Q$, such 
that $R$ is a left $T$-module and $\iota$ is left $T$-linear, then the
following statements hold. 
\begin{enumerate}[(i)]
\item $\Sigma':=R\ot_R\Sigma$ is an $R$-firmly projective right
  $A$-module; 
\item The functor $-\ot_R\Sigma\simeq -\ot_R\Sigma':\Mm_R\to\Mm^\Cc$ is fully
  faithful; 
\item If moreover the 
functor $\Hom^\Cc(\Sigma,-):\Mm^\Cc \to \Mm_T$ is fully faithful (as e.g. in
the setting of \cite[Theorem 4.1]{BohmVer:cleftex} or in forthcoming Theorem
\ref{thm:corext}), 
then $-\ot_R\Sigma:\Mm_R\to\Mm^\Cc$ is an equivalence. 
\end{enumerate}
\end{theorem}

\begin{proof}
\ul{(i)}. This assertion follows immediately by Corollary \ref{cor:firm_ideal}
(ii). 

\ul{(ii)}.
Applying \thref{Moritafirm} (ii) to the $R$-reduced Morita
  context $$(R,S:= \Sigma R\wt Q \subseteq {}^* \cC, R \otimes_R \Sigma, Q
  \otimes_R R, \wtb,\btb),$$ one concludes that the functor $F:= -\otimes_R
  \Sigma: \Mm_R \to \Mm_S$ is fully faithful. Moreover, $F$ factorizes as 
\[
\xymatrix{
\Mm_R \ar[rr]^-{ -\otimes_R \Sigma}&&
\Mm^\cC \ar[rr]&&
\Mm_{\cC^*}\ar[rr]&&
\Mm_S.}
\]
The second and third functors act on the morphisms as the identity map, hence
their composite is faithful, hence fully faithful. This proves that the
leftmost arrow describes a full, hence fully faithful functor.

\ul{(iii)}. By part (ii) we know that the functor $-\ot_R\Sigma:\Mm_R \to
\Mm^\Cc$ is fully faithful. By assumption, also 
$\Hom^\Cc(\Sigma,-):\Mm^\Cc\to\Mm_T$ is fully faithful. It was proven in
\cite{GomVer:ComatFirm} that this last statement is equivalent to the fact
that the functor $\Hom^\Cc(\Sigma,-)\ot_RR:\Mm^\Cc\to\Mm_R$ is fully faithful,
since $R$ is a left ideal in $T$.  
\end{proof}

\subsection{A Morita context associated to a pure coring extension}
\selabel{contextextension} 

In this section we consider two corings $\cD$ and $\Cc$ over unital rings $L$
and $A$, respectively, such that $\cC$ is a $\cC$-$\cD$ bicomodule via the
left regular $\cC$-coaction (i.e.\ $\cD$ is a {\em right extension} of $\cC$
in the sense of \cite{Brz:corext}). 
Assume that $\cD$ is a {\em pure} coring extension of $\cC$, in the sense
that, for any right $\cC$-comodule $(M,\varrho)$, the equalizer
$$
\xymatrix{
M \ar[rr]^-{\varrho}&&
M \ot_A {\mathcal C} \ar@<2pt>[rr]^-{\varrho\ot_A {\mathcal C}}
\ar@<-2pt>[rr]_-{M\ot_A \Delta_{\mathcal C}}&&
M \ot_A {\mathcal C} \ot_A {\mathcal C}
}
$$
in $\Mm_L$ is $\cD\ot_L \cD$-pure, i.e. it is preserved by the functor $(-)\ot_L
{\mathcal D}\ot_L {\mathcal D}:{\mathcal M}_L \to {\mathcal M}_L$.
Then, in addition to
\eqref{eq:A_mod_context} and 
\equref{contextSigma}, we can associate to $\Sigma\in {}_L \Mm^\cC$ a further
Morita context  
\begin{equation}\eqlabel{contextextension}
\widetilde{\MM}(\Sigma)=({}_L\Hom_L(\cD,T),{}^\cC\End^\cD(\cC)^{op},{}_L
\Hom^\cD(\cD,\Sigma),\widetilde{Q}, 
\wdi, \bdi ). 
\end{equation}
Here $T=\End^\cC(\Sigma)$, ${}_L\Hom_L(\cD,T)$ is a convolution algebra and 
${}^\cC\End^\cD(\cC)^{op}$ is the (opposite) endomorphism algebra of $\Cc$ as
$\Cc$-$\Dd$ bicomodule. The bimodule ${\widetilde Q}$ is a subset of $Q$ in
\equref{contextSigma}, for 
whose elements $q$ the right $L$-linearity  condition $q(x)(cl)=q(x)(c)l$
holds, for $x\in \Sigma$, $c\in \cC$ and $l\in L$. The connecting maps are
expressed in terms of the connecting maps in \equref{contextSigma} as
\begin{eqnarray} \label{eq:bd}
&&\\
&\bdi: 
{\widetilde Q}\ot_{{}_L\Hom_L(\cD,T)} {}_L\Hom^\cD(\cD,\Sigma) \to 
{}^\cC\End^\cD(\cC)^{op}, \quad &q
\otimes
p \mapsto \big(\ c\mapsto c_{[0]}(q\btd p(c_{[1]}))\ \big)\nonumber\\
&\wdi:{}_L\Hom^\cD(\cD,\Sigma)\ot_{{}^\cC\End^\cD(\cC)^{op}}
{\widetilde Q}\to {}_L\Hom_L(\cD,T), \quad 
&p
\otimes q \mapsto \big(\ d\mapsto p(d) \wtd
q\ \big), \nonumber
\end{eqnarray}
where a Sweedler type index notation $c\mapsto c_{[0]}\sstac L c_{[1]}$ is
used for the $\cD$-coaction in $\cC$ (implicit summation is understood). For
an explanation of the categorical 
origin of this Morita context and explicit form of the bimodule structures we
refer to \cite[Proposition 3.1]{BohmVer:cleftex} and its corrigendum. In order
to generalize in 
Theorem \ref{thm:corext} below some of the claims in Theorem 3.6 and Theorem
4.1 in \cite{BohmVer:cleftex} beyond the case when the connecting map $\bdi$
is surjective (in particular when $\Sigma$ is a {\em cleft} bicomodule), we
need the following  

\begin{lemma}\label{lem:com_firm}
Let $\cC$ be a coring over a unital ring $A$ and let $R\subseteq {}^*\cC$ be a
(non-unital) subring. 
If $\cC$ is a firm right $R$-module and the left regular $R$-module is flat,
then every right $\cC$-comodule $M$ is a firm right $R$-module, with action
$mr:=m^{[0]} r(m^{[1]})$.
\end{lemma}

\begin{proof}
For any $M\in \Mm^\cC$, there is a sequence of isomorphisms
\[
M\ot_R R \cong (M\ot^{\cC} \cC)\ot_R R \cong M\ot^\cC (\cC\ot_R R) 
\cong M\ot^\cC \cC \cong M
\]
mapping $M\ot_R R\ni m\ot_R r$ to $m^{[0]} r(m^{[1]})$. The second isomorphism
holds since $R$ is a flat left $R$-module and the penultimate isomorphism
holds since $\cC$ is a firm right $R$-module.
\end{proof}

The ring ${}^\cC\End^\cD(\cC)^{op}$ in the Morita context
\equref{contextextension} is a (unital) subring of the ring
${}^\cC\End(\cC)^{op}\cong {}^*\cC$.  Hence any right $\cC$-comodule $N$ is a
right ${}^\cC\End^\cD(\cC)^{op}$-module via 
\begin{equation}\label{eq:u_acts}
nu:= n^{[0]} \epsilon_\cC\big(u(n^{[1]})\big), \qquad \textrm{for } n\in N,\
u\in {}^\cC\End^\cD(\cC)^{op}.
\end{equation}
Obviously, the $\cC$-coaction $\rho^N:N\to N\ot_A\cC$ on $N$ is a right
${}^\cC\End^\cD(\cC)^{op}$-linear morphism, i.e.\ $(nu)^{[0]} \ot_A
(nu)^{[1]}=n^{[0]}\ot_A n^{[1]}u=n^{[0]}\ot_A n^{[1]}\varepsilon(u(n^{[2]}))$. 
 
\begin{theorem}\label{thm:corext}
Let $\cD$ be a coring over a unital ring $L$, which is a pure right extension
of a 
coring $\cC$ over a unital ring $A$. Let $\Sigma$ be an $L$-$\cC$ bicomodule
and consider the associated Morita context \equref{contextextension}. Let $R$
be a firm ring and a left ideal in ${\widetilde Q}\bdi
{}_L\Hom^\cD(\cD,\Sigma)$, such that the left regular $R$-module is flat and
$\cC$ is a firm right $R$-module. Then 
\begin{equation}\label{eq:can_nat}
\can: \Hom_A(\Sigma,-)\ot_T \Sigma \to -\ot_A \cC,\qquad 
\can_N(\phi_N \ot_T x) = \phi_N(x^{[0]})\ot_A x^{[1]}
\end{equation}
is a natural isomorphism 
and the functor $\Hom^\cC(\Sigma,-):\Mm^\cC\to \Mm_T$ is fully faithful. 
\end {theorem}

\begin{proof}
By Corollary \ref{cor:firm_ideal} there exists an $R$-bimodule map 
\[
R \to {\widetilde Q}\ot_{{}_L \Hom_L(\cD,T)}
{}_L\Hom^\cD(\cD,\Sigma), \qquad r\mapsto {\widetilde{\jmath}}_r \ot j_r
\]
(with implicit summation understood), such that ${\widetilde{\jmath}}_r \bdi
j_r=r$. 
We claim that the inverse of \eqref{eq:can_nat} is given by the well defined
map  
\[
\can_N^{-1}(n\ot_A cr)=n {\widetilde{\jmath}}_r(c_{[0]})(-)\ot_T j_r(c_{[1]}),
\qquad \textrm{for } n\in N,\ c\in \cC,\ r\in R.
\]
Indeed, the same arguments used to prove \cite[Theorem 3.6]{BohmVer:cleftex}
yield 
\begin{equation}\label{eq:can_inv}
(\can_N \circ \can_N^{-1})(n\ot_A cr)=n\ot_A cr\qquad \textrm{and} \qquad 
(\can_N^{-1} \circ \can_N) (\phi_N\ot_T xr) =\phi_N\ot_T xr,
\end{equation}
for $n\in N$, $c\in \cC$, $\phi_N \in \Hom_A(\Sigma,N)$, $x\in
\Sigma$ and $r\in R$,
where the $R$-actions are induced by \eqref{eq:u_acts}.
$\cC$ is a firm right $R$-module by assumption. $R$ is a non-unital subring in
${}^\Cc\End^\cD(\cC)^{op}\subseteq {}^\cC\End(\cC)^{op}\cong {}^*\cC$, hence
$\Sigma$ is a firm right $R$-module by Lemma \ref{lem:com_firm}. Thus
$N\ot_A \cC R=N\ot_A \cC$ and $\Hom_A(\Sigma,N)\ot_T \Sigma
R=\Hom_A(\Sigma,N)\ot_T \Sigma$. Therefore \eqref{eq:can_inv} proves that
\eqref{eq:can_nat} is a natural isomorphism.  

In view of \cite[Lemma 2.1 (2)]{BohmVer:cleftex}, for any right $\cC$-comodule
$N$ with coaction $\rho^N$ and $n\in N$, $r\in R$, $(\can_N^{-1} \circ
\rho^N)(nr)\in \Hom^\cC(\Sigma,N)\ot_T \Sigma$. Since $N$ is a firm
right $R$-module Lemma by \ref{lem:com_firm}, this shows that the range of
$\can_N^{-1} \circ \rho^N$ lies within $\Hom^\cC(\Sigma,N)\ot_T \Sigma$.
The same computations in \cite[Theorem 4.1]{BohmVer:cleftex} yield that
corestriction of $\can_N^{-1} \circ \rho^N$ gives the inverse of
\eqref{eq:T-counit},  hence $\Hom^\cC(\Sigma,-):\Mm^\cC\to \Mm_T$ is a fully
faithful functor.  
\end{proof}

\subsection{A Morita context connecting two comodules} \label{sec:com_Morita}
Two objects $\Sigma$ and $\Lambda$ in a k-linear category determine a Morita
context 
\begin{equation}\label{eq:Sigma-Lambda_context}
\mathbb{M}(\Sigma,\Lambda)=(\End(\Sigma), \End(\Lambda), \Hom(\Lambda,\Sigma),
\Hom(\Sigma,\Lambda), \square,\blacksquare),
\end{equation}
where multiplication, all bimodule structures and also the connecting maps are
given by 
composition in the category (what will be denoted by juxtaposition
throughout). In this section we study (reduction of) the Morita context
\eqref{eq:Sigma-Lambda_context}, determined by two objects $\Sigma$ and 
$\Lambda$ in the $k$-linear category of right comodules of a coring $\cC$ over
a unital ring $A$ over a commutative ring $k$.  
Throughout the section let $B\subseteq \Hom^\cC(\Sigma,\Lambda)
\Hom^\cC(\Lambda,\Sigma)$ be a left ideal and an idempotent ring
and $W:=\Hom^\cC(\Lambda,\Sigma) B \Hom^\cC (\Sigma, \Lambda)$.

We can consider the reduced form \eqref{eq:red_context} (or equivalently, 
\equref{2ndreduced}, see Remark \ref{rem:B-W_eq} (iii)) of the Morita context
\eqref{eq:Sigma-Lambda_context}, i.e.\ 
\begin{equation}\eqlabel{SigmaLambdareduced}
(W,B,\Hom^\cC(\Lambda,\Sigma)\ot_B B, B\ot_B \Hom^\cC(\Sigma,\Lambda),
  \ol{\square},\ol{\blacksquare}). 
\end{equation}
We obtain the following (not necessarily commutative) diagram of adjoint
functors. (In order to see that $G_\Sigma$ and $G_\Lambda$ are well defined,
consult \thref{idempotentfirm} (i).)
\begin{equation}\label{eq:triangle_diagram}
\xymatrix{
&&&& \Mm^\Cc \ar@<-.5ex>[ddllll]_-{G_\Sigma=\Hom^\Cc(\Sigma,-)W\ot_WW
    \hspace{1cm}}  
\ar@<.5ex>[ddrrrr]^-{\hspace{1cm} G_\Lambda=\Hom^\Cc(\Lambda,-)B\ot_BB} \\
{}\\
\Mm_{\widetilde{W}}\cong\Mm_W \ar@<-.5ex>[uurrrr]_-{F_\Sigma=-\ot_W\Sigma}
\ar@<.5ex>[rrrrrrrr]^-{-\ot_W\Hom^\Cc(\Lambda,\Sigma)\ot_BB} 
&&&&&&&& \Mm_B\cong\Mm_{\widetilde{B}}
\ar@<.5ex>[uullll]^-{F_\Lambda=-\ot_B\Lambda}
\ar@<.5ex>[llllllll]^-{-\ot_B\Hom^\Cc(\Sigma,\Lambda)} \ ,
}
\end{equation}
where ${\widetilde{W}}= W\ot_W W$, as before.
If $W$ is a firm ring (i.e.\ $W={\widetilde W}$) then the adjunction
$(F_\Sigma,G_\Sigma)$ reduces to  \equref{HomTensorSigma}.
If we consider $F_\Sigma$ and  $G_\Sigma$ as functors between $\Mm^\Cc$ and
$\Mm_W$, the unit and counit are defined as
\begin{eqnarray*}
&\nu^\Sigma_N:N \to \Hom^\cC(\Sigma,N\ot_W \Sigma)W\ot_W W,\quad\ \quad 
&n\mapsto ((n^w)^{w'} \ot_W -)w'\ot_W w,\\
&\zeta^\Sigma_M: \Hom^\cC(\Sigma,M)W\ot_W W \ot_W \Sigma\to M,\qquad 
&\phi w \ot_W w'\ot_W x \mapsto \phi(ww'x),
\end{eqnarray*}
for $M\in \Mm^\cC$ and $N\in \Mm_W$. Similarly, we define the unit
$\nu^\Lambda$ and the counit $\zeta^\Lambda$ for the adjunction
$(F_\Lambda,G_\Lambda)$. 

The aim of Proposition \ref{prop:nat_eq} is to relate the functors
$G_\Sigma:\Mm^\cC \to \Mm_W$ and $G_\Lambda:\Mm^\cC \to \Mm_B$, i.e.\ to show
that the outer triangle in diagram \eqref{eq:triangle_diagram} is commutative
up to a natural isomorphism.

\begin{proposition} \label{prop:nat_eq}
For a unital ring $A$, 
let $\Sigma$ and $\Lambda$ be right comodules of an $A$-coring $\cC$ and
$B \subseteq \End^\cC(\Lambda)$ and $W\subseteq \End^\cC(\Sigma)$ as
above. Then, for any right $\cC$-comodule $M$, there is a right $W$-module
isomorphism 
\[
\Hom^\cC(\Sigma,M)W\ot_W W \cong \Hom^\cC(\Lambda,M)B \ot_B B\ot_B
\Hom^\cC(\Sigma,\Lambda).
\]
\end{proposition}

\begin{proof}
Consider the $B$-reduced form \equref{SigmaLambdareduced}
of the Morita context \eqref{eq:Sigma-Lambda_context}. 
The right $W$-module $\Hom^\cC(\Sigma,M)W\ot_W W$ and the right $B$-module 
$\Hom^\cC(\Lambda,M)B\ot_B B$ are firm
by \thref{idempotentfirm} (i), for any right $\cC$-comodule $M$. Therefore,
by \leref{omegabeta}, the morphisms  
$\Hom^\cC(\Sigma,M)W\ot_W W\ot_W \ol{\square}$ and
$\Hom^\cC(\Lambda,M)B\ot_B B \ot_B \ol{\blacksquare}$ are
isomorphisms. Furthermore, composition of 
$\cC$-comodule morphisms defines maps, for any $M\in \Mm^\cC$,
\begin{eqnarray*}
&&\omega_1: \Hom^\cC(\Sigma,M)W\ot_W W \ot_W \Hom^\cC(\Lambda,\Sigma)\ot_
B B \to \Hom^\cC(\Lambda,M)B,\\
&&\omega_2: \Hom^\cC(\Lambda,M)B\ot_B B \ot_B B \ot_
B\Hom^\cC(\Sigma,\Lambda) \to \Hom^\cC(\Lambda,M)B\Hom^\cC(\Sigma,\Lambda),\\
&&\omega_3: \Hom^\cC(\Lambda,\Sigma)\ot_B B \ot_B \Hom^\cC(\Sigma,\Lambda)
\to W.
\end{eqnarray*}
Obviously, $\omega_1$ is a right $B$-module map, $\omega_2$ is right
$W$-linear and $\omega_3$ is $W$-$W$ bilinear. By Lemma \ref{lem:B-W_prop}
(ii), $\Hom^\cC(\Lambda,M)B\Hom^\cC(\Sigma,\Lambda)=\Hom^\cC(\Lambda,M)B
\Hom^\cC(\Sigma,\Lambda)W$ is a right $W$-submodule of
$\Hom^\cC(\Sigma,M)W$. Hence there is a well defined map  
\begin{eqnarray*}
(\omega_2 \ot_W \omega_3)&\circ&\big((\Hom^\cC(\Lambda,M)B  \ot_
B B \ot_B \ol{\blacksquare})^{-1} \ot_B
\Hom^\cC(\Sigma,\Lambda)\big)\\ 
&\circ& \big((\Hom^\cC(\Lambda,M)B  \ot_B \mu_B)^{-1}\ot_B
\Hom^\cC(\Sigma,\Lambda)\big) 
\end{eqnarray*}
from $\Hom^\cC(\Lambda,M)B \ot_B B\ot_B \Hom^\cC(\Sigma,\Lambda)$ to
$\Hom^\cC(\Sigma,M)W \ot_W W$. 
A routine computation shows that it is an isomorphism with inverse
\[
(\omega_1\ot_B B\ot_B \Hom^\cC(\Sigma,\Lambda))\circ 
(\Hom^\cC(\Sigma,M)W \ot_W W\ot_W \ol{\square})^{-1}\circ
(\Hom^\cC(\Sigma,M)W \ot_W \mu_W)^{-1}.
\]
This ends the proof.
\end{proof}

\begin{corollary}\label{cor:F_Lambda-F_Sigma}
Let $\Sigma$ and $\Lambda$ be right comodules of a coring $\cC$ over a unital
ring $A$, and let
$B \subseteq \End^\cC(\Lambda)$ and $W\subseteq \End^\cC(\Sigma)$ be non-unital
subrings as in Proposition \ref{prop:nat_eq}. Then the following assertions
hold. 
\begin{enumerate}[(i)]
\item If $\zeta_\Sigma^\Lambda$ or $\zeta_\Lambda^\Sigma$ is an isomorphism,
  then the functor $F_\Sigma:\Mm_W\to\Mm^\Cc$ is fully faithful if and only if
  $F_\Lambda:\Mm_B\to\Mm^\Cc$ is fully faithful;
\item The functor $G_\Sigma:\Mm^\cC \to \Mm_W$ is fully faithful if and
only if $G_\Lambda:\Mm^\cC \to \Mm_B$ is fully faithful;
\item The functor $G_\Sigma:\Mm^\cC \to \Mm_W$ is an equivalence if
  and only if $G_\Lambda:\Mm^\cC \to \Mm_B$ is an equivalence.
\end{enumerate}
\end{corollary}

\begin{proof}
\ul{(i)}. For any $N\in\Mm_W$, there is a natural morphism
$$
N\ot_W\zeta_\Sigma^\Lambda:N\ot_W\Hom^\Cc(\Lambda,\Sigma)\ot_B B \ot_B
\Lambda\to N\ot_W\Sigma.
$$ 
Therefore, $F_\Sigma$ is naturally isomorphic to the composite of the functors
$F_\Lambda$ and $-\ot_W\Hom^\Cc(\Lambda,\Sigma)\ot_BB$, provided
$\zeta^\Lambda_\Sigma$ is an isomorphism. Since we know that
$-\ot_W\Hom^\Cc(\Lambda,\Sigma)\ot_BB:\Mm_W\to\Mm_B$ is an equivalence (see
Remark \ref{rem:B-W_eq} applied to the $B$-reduced form of the Morita context 
\eqref{eq:Sigma-Lambda_context}), this proves the claim.\\
\ul{(ii)\& (iii)}.
By Proposition \ref{prop:nat_eq}, $G_\Sigma$ is naturally isomorphic to the
composite of $G_\Lambda$ and the equivalence functor $-\ot_B
\Hom^\cC(\Sigma,\Lambda):\Mm_B \to \Mm_W$,
which proves both claims. 
\end{proof}

In \cite[Proposition 2.7]{BohmVer:cleftex} we proved that, in the case when
$\cC$ is a finitely generated projective left $A$-module, the Morita context
$\mathbb{M}(\Sigma)$ in \equref{contextSigma} is strict if and only if the
Strong Structure Theorem holds, that is, $\Hom^\cC(\Sigma,-): \Mm^\cC\to
\Mm_{T}$ is an equivalence. The aim of the rest of current section is to
extend this result beyond the case when $\cC$ is a finitely generated
projective left $A$-module. 
    
In order to apply the results of this section,
in addition to
$\Sigma$ we need a second $\cC$-comodule. In what follows we give sufficient
and necessary conditions under which 
the range $B$ of the connecting map $\bt$ in the Morita context
\equref{contextSigma} has a $B$-$\cC$ bicomodule structure such that 
the corresponding adjunction $(F_B,G_B)$ (see \eqref{eq:triangle_diagram}) is
an equivalence. In the case when these conditions hold, we apply Corollary
\ref{cor:F_Lambda-F_Sigma} to prove that also the adjunction
$(F_\Sigma,G_\Sigma)$ is an equivalence. 

Recall (e.g. from \cite{Wis:com_cat} or \cite{Ver:PhD}) that for a left module
$P$ over a unital ring $A$, the {\em finite topology} on
${}^*P:={}_A\Hom(P,A)$ is generated by the open sets ${\mathcal
  O}(f,p_1,\dots,p_n)=\{\ g\in {}^* P\ |\ g(p_i)=f(p_i),\ i=1,\dots,n\
\}$. The left $A$-module $P$ is said to be {\em   weakly locally projective}
if every finitely generated submodule of $P$ has a dual basis in $P\times {}^*
P$. Equivalently, if and only if ${}^*P$ satisfies the {\em
  $\alpha$-condition}, meaning that the map 
\[
M\ot_A P \to \Hom_{Ab}({}^* P,M),\qquad m\ot_A p \mapsto \big(\ f\mapsto m
f(p)\ \big)
\]
is injective, for every right $A$-module $M$. A non-unital ring $B$ has right
{\em local units} if for any finite subset $\{b_1,\dots,b_n\}$ of $B$
there exists an element $e\in B$ such that $b_ie=b_i$, for all
$i=1,\dots,n$. If $B$ is a ring with right local units then it is in
particular firm and its left regular module is flat.

\begin{lemma}\label{lem:G_B_tilde}
For a unital ring $A$, 
let $\Sigma$ be a right comodule of an $A$-coring $\cC$ and let $B:=Q\bt
\Sigma$ be the range of the connecting map $\bt$ in the Morita context
$\MM(\Sigma)$ in \equref{contextSigma}. 
Assume that the left regular $B$-module extends to a $B$-$\cC$ bicomodule such
that the connecting map $\bt$ corestricts to a $B$-$\cC$ bicomodule
epimorphism $Q\ot_{T} \Sigma\to B$. Then $F_B$ has a left
inverse, the functor $\widetilde{F}_B$, sending a right $\cC$-comodule $M$ to
the right $B$-module $M$, with action $mb:=m^{[0]} b(m^{[1]})$, and acting on
the morphisms as the identity map.
\end{lemma}

\begin{proof}
For $q\in Q$, $x\in \Sigma$, $b\in B$ and $c\in \cC$,
\begin{eqnarray*}
\big(q(x)^{[0]} b(q(x)^{[1]})\big) (c)&=&
q(x^{[0]})(c) b(x^{[1]})=
b\big(q(x^{[0]})(c)x^{[1]}\big)\\
 &=& b\big(c^{(1)} q(x)(c^{(2)})\big) =
\big(q(x) b)(c).
\end{eqnarray*}
The first equality follows by the right $\cC$-colinearity of
$\bt:q\ot x \mapsto q(x)$ and the form of the right $A$-action on ${}^*\cC$. 
The second equality follows by the left $A$-linearity of $b\in B\subseteq
{}^*\cC$. The penultimate equality follows by the defining property of $q\in
Q$ while the last one follows by the form of the multiplication in
$^*\cC$. Since $B$ is the range of $\bt$, we conclude that $b^{[0]}
b'(b^{[1]})=bb'$, for all $b,b'\in B$. Thus $\widetilde{F}_B\circ F_B$ takes a
firm right $B$-module $N$ to the right $B$-module $N$, with action
\[
n\ot b' \mapsto n^b b^{[0]} b'(b^{[1]}) = n^b b b' =n b',
\]
where $\sd_{N,B}(n)= n^b\ot_B b$ is the unique element of $N\ot_B B$ such that
$n^bb=n$. This proves $\widetilde{F}_B\circ F_B=\Mm_B$. 
\end{proof}

\begin{theorem}\label{thm:Morita_G_B_eq}
For a unital ring $A$, 
let $\Sigma$ be a right comodule of an $A$-coring $\cC$ and let $B:=Q\bt
\Sigma$ be the range of the connecting map
$\bt$ in the Morita context $\MM(\Sigma)$ in \equref{contextSigma}. The
following assertions are equivalent.
\begin{enumerate}[(i)]
\item The left regular $B$-module extends to a $B$-$\cC$ bicomodule such that
  the connecting map $\bt$ corestricts to a $B$-$\cC$ bicomodule epimorphism
  $Q\ot_{T} \Sigma\to B$, and $F_B: \Mm_B\to \Mm^\Cc$ is an
  isomorphism;  
\item $\cC$ is weakly locally projective as a left $A$-module and $B$
  is dense in the finite topology on ${}^*\cC$;
\item $\cC$ is weakly locally projective as a left $A$-module, $B$ has
  right local units (in particular, $B$ is a firm ring) and $\cC$ is firm as a
  right $B$-module;
\item The left regular $B$-module extends to a $B$-$\cC$ bicomodule such that
  the connecting map $\bt$ corestricts to a $B$-$\cC$ bicomodule epimorphism
  $Q\ot_{T} \Sigma\to B$, 
  $\Cc$ is firm as a right $B$-module and the left regular $B$-module is flat.
\end{enumerate} 
\end{theorem}

\begin{proof} 
\ul{(i)$\Rightarrow$ (ii)}. 
By Lemma \ref{lem:G_B_tilde}, 
$F_B$ has a left inverse, the functor $\widetilde{F}_B$, sending a right
$\cC$-comodule $M$ to the right $B$-module $M$, with action $mb:=m^{[0]}
b(m^{[1]})$, and acting on the morphisms as the identity map. 
Under assumptions (i), $F_B$ is an isomorphism. Hence
$\widetilde{F}_B=F_B^{-1}$.  
The category of firm modules over an idempotent ring was proven to be a
Grothendieck category by Mar\'{\i}n in \cite{Mar:Morita}. Since
$F_B^{-1}(\cC)=\cC$ is a firm right $B$-module (with $B$-action $cb=c^{(1)} 
b(c^{(2)})$), the smallest Grothendieck subcategory $\sigma[\cC_B]$ of
$\widetilde{\Mm}_B$, which contains $\cC$, is contained in $\Mm_B$. On the
other hand, by \cite[Corollary 3.30]{Ver:PhD}, any right $\cC$-comodule is
subgenerated by $\cC$ as a right $B$-module, i.e.\ $\Mm^\cC$ is contained in
$\sigma[\cC_B]$. Hence the isomorphism 
$\Mm^\cC\cong \Mm_B$ implies $\sigma[\cC_B]\cong \Mm^\cC$. 
One can easily adapt the proof of \cite[Theorem 3.5
  (a)$\Rightarrow$(d)]{Wis:com_cat} to conclude that the map 
\[
\alpha_{P,B}:P\ot_A \cC \to \Hom_{Ab}(B,P),\qquad p\ot_A c \mapsto \big(\
b \mapsto pb(c)\ \big)
\]
is injective, for any right $A$-module $P$. This is equivalent to assertion
(ii) by \cite[Theorem 2.58]{Ver:PhD}. 

\ul{(ii)$\Leftrightarrow$(iii)}. This equivalence is proven in \cite[Corollary
2.48]{Ver:PhD}. 

\ul{(iii)$\Rightarrow$(iv)}. Since $B$ has right local units, its left
regular module is flat. 
The existence of the required $B$-$\cC$ bicomodule
structure on $B$ follows by a rationality argument. 
By construction, the map $\bt:Q\ot_T\Sigma\to B$ is a surjective $B$-$\*C$
bilinear map. 
Since $Q\ot_T\Sigma$ is a $B$-$\Cc$ bicomodule, it is in particular a rational
right $\*C$-module. Hence, $B$ being the image of the map $\bt$, it is a
quotient 
of the rational $\*C$-module $Q\ot_T\Sigma$, and hence $B$ itself is rational
by \cite[Proposition 4.2]{CVW:rational}. Therefore $B$ is a $B$-$\Cc$
bicomodule and $\bt$ corestricts to a $B$-$\Cc$ bicomodule map. 

\noindent

\ul{(iv)$\Rightarrow$(i)}. By Lemma \ref{lem:G_B_tilde}, $F_B$ has a left
inverse $\widetilde{F}_B$. Composition $F_B\circ \widetilde{F}_B$ makes sense
by 
Lemma \ref{lem:com_firm}.
The proof is completed by computing the coaction on the right $\cC$-comodule
$F_B\circ \widetilde{F}_B(M)$, for $M\in \Mm^\cC$. 
For $q\in Q$, $x\in \Sigma$ and $c\in \cC$,
\[
q(x)^{[0]} (c) q(x)^{[1]} = q(x^{[0]})(c) x^{[1]} = c^{(1)} q(x)(c^{(2)}).
\]
The first equality follows by right $\cC$-colinearity of $\bt$, and the
second equality follows by the defining property of $q\in Q$. Hence, for
$b\in B$ and $c\in \cC$, $b^{[0]} (c) b^{[1]} = c^{(1)}
b(c^{(2)})$. Note that by Lemma \ref{lem:com_firm} $\widetilde{F}_B(M)=
\widetilde{F}_B(M)B$. The right $\cC$-coaction on $F_B\circ
\widetilde{F}_B(M)$ comes out as 
\begin{eqnarray*}
mb&\mapsto& m^{[0]} b^{[0]}(m^{[1]})\ot_A b^{[1]}=
m^{[0]} \ot_A b^{[0]}(m^{[1]}) b^{[1]}=
m^{[0]} \ot_A m^{[1]} b(m^{[2]}) =\\
&&(m^{[0]} b(m^{[1]}))^{[0]} \ot_A (m^{[0]} b(m^{[1]}))^{[1]}=
(m b)^{[0]}\ot_A (m b)^{[1]}.
\end{eqnarray*}
This proves $F_B\circ \widetilde{F}_B=\Mm^\cC$ hence the theorem.
\end{proof}

For a unital ring $A$, let $\Sigma$ be a right comodule of an $A$-coring $\cC$
and let $B:=Q\bt 
\Sigma$ be the range of the connecting map $\bt$ in the Morita context
$\MM(\Sigma)$ in \equref{contextSigma}. Assume that the equivalent
conditions in Theorem \ref{thm:Morita_G_B_eq} hold. Then $B$ is a firm ring,
and we can consider the $B$-reduced form of $\MM(\Sigma)$
\begin{equation}\eqlabel{MSigmareduced}
\big((\Sigma\wt Q)(\Sigma\wt Q)\ ,\  B\ ,\ \Sigma \ot_B B\ ,\ B \ot_
  B Q\ ,\ \widetilde{\wt}\ ,\ \widetilde{\bt}\big ) ,
\end{equation}
where the connecting maps are given, for $b,\widetilde{b},\in B$, $x\in
\Sigma$, $q\in Q$, by 
\begin{eqnarray}
&(x \ot_B \widetilde{b}) \widetilde{\wt} (b\ot_B q) :=
  x \widetilde{b}  \wt b q \qquad \textrm{and}\qquad 
&(b\ot_B q) \widetilde{\bt} (x \ot_B \widetilde{b}) :=
  b q \bt x \widetilde{b}. 
\label{eq:conn_maps1}
\end{eqnarray}

On the other hand, under the conditions in Theorem \ref{thm:Morita_G_B_eq}, $B$
is also a right $\cC$-comodule, hence we can consider a Morita 
context $\mathbb{M}(\Sigma,B)$ as in \eqref{eq:Sigma-Lambda_context}. 
In the next lemma we show that 
also the Morita context $\mathbb{M}(\Sigma,B)$ admits a $B$-reduced form.

\begin{lemma}\label{lem:B_left_ideal}
For a unital ring $A$, let $\Sigma$ be a right comodule of an $A$-coring $\cC$
and let $B:=Q\bt \Sigma$ be the range of the connecting map $\bt$ in the
Morita context $\MM(\Sigma)$ in \equref{contextSigma}. Assume that the
equivalent conditions in Theorem \ref{thm:Morita_G_B_eq} hold, hence $B$ is a
$B$-$\cC$ bicomodule.
Then $B$ is a left ideal in $\Hom^\Cc(\Sigma,B)\Hom^\Cc(B,\Sigma)$.   
\end{lemma}

\begin{proof}
Note first that there is a well-defined map
\begin{equation}\label{eq:gamma_def}
\gamma:\Sigma\to\Hom^\Cc(B,\Sigma),\qquad \gamma(y)(b)=yb \ .
\end{equation}
That is, for all $y\in\Sigma$, the map $\gamma(y):B\to \Sigma,\
\gamma(y)(b)=yb=y^{[0]} b(y^{[1]})$ is right $\Cc$-colinear. Indeed, for $q\in
Q$ and  $x\in \Sigma$, 
\begin{eqnarray*}
\gamma(y)(q(x)^{[0]})\ot_A q(x)^{[1]} &=&
y q(x)^{[0]} \ot_A q(x)^{[1]} =
y^{[0]}\ot_A q(x)^{[0]}(y ^{[1]}) q(x)^{[1]} \\
&=&y^{[0]} \ot_A q(x^{[0]})(y ^{[1]}) x^{[1]} =
y^{[0]} \ot_A y^{[1]} q(x)(y^{[2]}) \\
&=&\big(y^{[0]} q(x)(y^{[1]})\big) ^{[0]} \ot_A \big(y^{[0]}
q(x)(y^{[1]})\big) ^{[1]}\\
 &=&\gamma(y)(q(x))^{[0]}\ot_A \gamma(y)(q(x))^{[1]}.
\end{eqnarray*}
The second and the last equalities follow by the form of the $B$-action on
$\Sigma$. The 
third equality follows by the right $\cC$-colinearity of $\bt$ and the fourth
equality is a consequence of the defining property of $q\in Q$. The
penultimate equality follows by right $A$-linearity of the $\cC$-coaction on 
$\Sigma$.

Next, since any morphism in $\End^\Cc(B)$ is right $B$-linear,
the map $\beta: B\to\End^\Cc(B)$, $\beta(b)(b')=bb'$ turns $B$ into a
left ideal in $\End^\Cc(B)$. Remark furthermore that
$\Hom^\Cc(\Sigma,B)\Hom^\Cc(B,\Sigma)$ is in a natural way a (two-sided) ideal
in $\End^\Cc(B)$. Finally, for any $q\bt x\in B$, by right
  $B$-linearity of $q\in Q$, we have that 
$\beta(q\bt x)=q\gamma(x)\in\Hom^\Cc(\Sigma,B)\Hom^\Cc(B,\Sigma)$. Hence the
image of $\beta$ and therefore $B$ is a left ideal in
$\Hom^\Cc(\Sigma,B)\Hom^\Cc(B,\Sigma)$. 
\end{proof}

In view of Lemma \ref{lem:B_left_ideal}, we can apply the theory developed 
at the beginning of this section to obtain a $B$-reduced form of the Morita
context $\MM(\Sigma,B)$, as in
\equref{SigmaLambdareduced}. Note that since $B$ is a firm ring, $M\ot_ B B
\cong MB\ot_B B$ and $B\ot_B N \cong B\ot_B BN$, for any $M\in
      {\widetilde{\Mm}}_B$ and $N\in {}_B {\widetilde{\Mm}}$. Thus from
      \equref{SigmaLambdareduced} we obtain
\begin{equation}\eqlabel{BSigmareduced}
\big(\Hom^\cC(B,\Sigma)B \Hom^\cC(\Sigma,B) ,\ B  ,\ \Hom^\cC(B,\Sigma) \ot_
  B B  ,\  B\ot_B  \Hom^\cC(\Sigma,B),\  \widetilde{\wdi},\
  \widetilde{\bdi}  \big) ,
\end{equation}
where the connecting maps are given, for $b,\widetilde{b}\in B$, $\zeta\in
\Hom^\cC(\Sigma,B)$ and $\xi \in \Hom^\cC(B,\Sigma)$, by
\begin{eqnarray}
&(\xi \ot_B b) \,\widetilde{\wdi}\, (\widetilde{b} \ot_B \zeta) := \xi
  b\widetilde{b}\zeta \quad\qquad 
  \textrm{and}\qquad  
&(\widetilde{b} \ot_B  \zeta)\, \widetilde{\bdi}\, (\xi \ot_B
  b) := \widetilde{b}\zeta \xi b.
\label{eq:conn_maps}
\end{eqnarray}

\begin{proposition}\label{prop:two_contexts}
For a unital ring $A$, let $\Sigma$ be a right comodule of an $A$-coring $\cC$
and let $B:=Q\bt 
\Sigma$ be the range of the connecting map $\bt$ in the Morita context
$\MM(\Sigma)$ in \equref{contextSigma}. Assume that the equivalent
conditions in Theorem \ref{thm:Morita_G_B_eq} hold, hence $B$ is a $B$-$\cC$ 
bicomodule. Then the Morita contexts \equref{MSigmareduced} and
\equref{BSigmareduced} are isomorphic. 
\end{proposition}

\begin{proof}
In terms of the map $\gamma$ in \eqref{eq:gamma_def}, put
\[
\alpha:\Sigma \ot_BB \to \Hom^\cC(B,\Sigma) \ot_BB,\qquad 
x\ot_B b\mapsto \gamma(x)\ot_B b.
\]
It is an isomorphism with inverse
\[
\alpha^{-1}: \Hom^\cC(B,\Sigma) \ot_B B \to \Sigma \ot_BB,\qquad 
\xi \ot_B bb' \mapsto \xi(b)\ot_B b'.
\]
Indeed, 
\begin{eqnarray*}
(\alpha^{-1}\circ\alpha)(x\ot_B bb')&=&\gamma(x)(b)\ot_B b'=xb\ot_B b'=
  x\ot_Bbb;\\  
(\alpha\circ\alpha^{-1})(\xi\ot_B bb')&=&\gamma(\xi(b))\ot_Bb' = \xi
  b\ot_B b' = \xi \ot_B bb',
\end{eqnarray*}
where in the penultimate equality of the second computation we used that a
right $\Cc$-comodule map $\xi$ is a right module map for $B\subseteq
{}^*\Cc$, hence
for all $b'\in B$, $\gamma(\xi(b))(b')=\xi(b)b'=\xi(bb')=(\xi b)(b')$.

By right $\cC$-colinearity of $\bt$, $Q\subseteq
\Hom^\cC(\Sigma,B)$. Conversely, for any $\zeta\in \Hom^\cC(\Sigma,B)$ and
$y\in \Sigma$, there exist (non-unique) elements $q_i\in Q$ and $x_i\in
\Sigma$ such that $B\ni\zeta(y) =\sum_i q_i(x_i)$. Thus, for $c\in \cC$, 
\begin{eqnarray*}
\zeta(y^{[0]})(c) y^{[1]} &=&
\zeta(y)^{[0]}(c)\zeta(y)^{[1]}
= \sum_i q_i(x_i)^{[0]}(c)q_i(x_i)^{[1]}\\
&=& \sum_i q_i(x_i^{[0]})(c) x_i^{[1]}
= \sum_i c^{(1)} q_i(x_i)(c^{(2)})
= c^{(1)} \zeta(y)(c^{(2)}).
\end{eqnarray*}
The first equality follows by the $\cC$-colinearity of 
$\zeta\in \Hom^\cC(\Sigma,B)$.
The third equality follows by the $\cC$-colinearity of $q_i$, for all
values of the index $i$. In order to conclude the penultimate equality we used
the defining property of $q_i\in Q$, for any index $i$. This proves that 
$\Hom^\cC(\Sigma,B) \subseteq Q$, hence the obvious map
\[
\beta:B\ot_B Q \to B\ot_B \Hom^\cC(\Sigma,B), \qquad 
b\ot_B q \mapsto b\ot_B q
\]
establishes an isomorphism. One checks easily that the isomorphisms $\alpha$
and $\beta$ are compatible with the connecting maps \eqref{eq:conn_maps1} and
\eqref{eq:conn_maps}. Thus in particular the ranges of the connecting maps
\eqref{eq:conn_maps1} and \eqref{eq:conn_maps} are coinciding (non-unital)
subrings of the endomorphism rings $\End^\cC(B)$ and $T$, respectively. That
is, 
$(\Sigma\wt Q) (\Sigma\wt Q)= \Hom^\cC(B,\Sigma) B \Hom^\cC(\Sigma,B)$. The
proof is completed by checking 
the bimodule map properties of $\alpha$ and $\beta$, what is
left to the reader. 
\end{proof}

The following theorem generalizes \cite[Theorem 5.22]{Ver:PhD} and hence
\cite[Theorem 5.3]{CVW:rational} and \cite[Theorem 4.15]{CaDeGrVe:comcor}
beyond the case when $\wt$ is surjective. 

\begin{theorem}\label{thm:str_str_thm}
For a unital ring $A$, let $\Sigma$ be a right comodule of an $A$-coring $\cC$
and consider the 
Morita context $\MM(\Sigma)$ in \equref{contextSigma}. Put $B:=Q\bt
\Sigma$ and $W:=(\Sigma\wt Q)(\Sigma\wt Q)$. Assume that the equivalent
conditions in Theorem \ref{thm:Morita_G_B_eq} hold. Then 
\begin{equation}\label{eq:F_Sigma}
G_\Sigma:= \Hom^\cC(\Sigma, -)W\ot_W W:\Mm^\cC\to \Mm_W 
\end{equation}
is an equivalence.
\end{theorem}

\begin{proof}
Since $F_B$ is an equivalence functor by assumption, so is its adjoint
$G_B$. 
By Lemma \ref{lem:B_left_ideal} $B$ is a
left ideal in $\Hom^\Cc(\Sigma,B)\Hom^\Cc(B,\Sigma)$
and by assumption $B$ is a firm ring. 
Thus the claim is an immediate consequence of Proposition
\ref{prop:two_contexts} and Corollary  \ref{cor:F_Lambda-F_Sigma}.
\end{proof}

If the functor \eqref{eq:F_Sigma} is fully faithful (as e.g. in Theorem
\ref{thm:str_str_thm}) then $\Sigma$ is a generator in $\Mm^\Cc$. 
Under the conditions in Theorem \ref{thm:str_str_thm}, $\Cc$ is a flat left
$A$-module (hence the forgetful functor $\Mm^\Cc\to \Mm_A$ preserves and
reflects monomorphisms) and \eqref{eq:F_Sigma} (thus also $F_\Sigma=-\ot_W
\Sigma:\Mm_W \to \Mm^\Cc$) is an equivalence. Hence $\Sigma$ is a faithfully
flat left $W$-module. 
Moreover, the following proposition holds. 

\begin{proposition}\label{prop:firm_proj}
In the situation described in Theorem \ref{thm:str_str_thm}, the functor
$-\ot_W \Sigma:\Mm_W \to \Mm_A$ has a right adjoint, the functor
\begin{equation}\label{eq:right_ad}
\Hom_A(\Sigma,-)W \ot_W W \simeq  -\ot_A \Hom_A(\Sigma,A)W \ot_W W.
\end{equation}
\end{proposition}

Note that if $W$ is a unital ring then the equivalence of the two forms of the
functor in \eqref{eq:right_ad} is equivalent to finitely generated
projectivity of $\Sigma$ as a right $A$-module. If $W$ is a firm ring and
$\Sigma$ is a firm left $W$-module, then this equivalence is equivalent to
$W$-firm projectivity of the right $A$-module $\Sigma$ by Theorem
\ref{thm:firm_proj}. 

\begin{proof}[Proof of Proposition \ref{prop:firm_proj}]
The functor $\Hom_A(\Sigma,-)W \ot_W W$ is equal to the composite of $-\ot_A
\cC:\Mm_A \to \Mm^\cC$ and the equivalence functor $G_\Sigma:\Mm^\cC\to \Mm_W$
in \eqref{eq:F_Sigma}. Since both of these functors possess as well a left
adjoint as a right adjoint, also  $\Hom_A(\Sigma,-)W \ot_W W$ possesses both
left and right adjoints. 
Furthermore, by \thref{idempotentfirm} (vi), there are natural equivalences
(of left adjoint functors) 
\begin{eqnarray*}
&&\Hom_A(\Sigma,-)W \ot_W W \simeq \Hom_A(\Sigma,-)
  \ot_{\widetilde{W}}{\widetilde{W}  }
\qquad\qquad\qquad\qquad \textrm{and}\quad\\
&&-\ot_A \Hom_A(\Sigma,A)W \ot_W W\cong -\ot_A \Hom_A(\Sigma,A) \ot_
{\widetilde{W}}{ \widetilde{W}}, 
\end{eqnarray*}
where $\widetilde{W}=W\ot_W W$ as before.
By \cite[Theorem 3.1 (i)$\Rightarrow$ (iii)]{Ver:equi} there is a natural
equivalence 
\[
\Hom_A(\Sigma,-) \ot_{\widetilde{W}  } {\widetilde{W}  }\cong -\ot_A
\Hom_A(\Sigma,A) \ot_{\widetilde{W}  } {\widetilde{W}  }.
\] 
We conclude the claim by combining these isomorphisms. 
\end{proof}

\begin{example}
Let $\Cc$ be an $A$-coring that is locally projective as left $A$-module and
put $B=\Rat(\*C)$, the rational part of the left dual of
$\Cc$. In several situations $B$ is dense is the finite
topology of $\Cc$. E.g. if $\Cc$ is a locally Frobenius coring as defined in
\cite{IV:cofrob}. If the base ring $A$ is a PF-ring, then the definition of
locally Frobenius coring is equivalent to the definition of a co-Frobenius
coring. If $A$ is a QF ring, then $B$ is dense in the finite
topology of $\*C$ if and only if $\cC$ is a semiperfect coring
  \cite{CI:semiperfect}.   

Morita contexts of type \equref{contextSigma} s.t. $\bt$ is surjective onto
$B=\Rat(\*C)$ have been considered extensively in
e.g. \cite{CaDeGrVe:comcor}, \cite{CVW:rational}, \cite{BDR} and fit into the
framework of this section. 
\end{example}

\subsection{Coseparable corings}\selabel{cosep}

Recall that an $A$-coring $\Cc$ is said to be \emph{coseparable} if and only
if there exists a $\Cc$-bicolinear left inverse $\mu:\Cc\ot_A\Cc\to \Cc$ 
of the comultiplication $\Delta$. If we denote 
$\gamma:=\varepsilon\circ\mu:\Cc\ot_A\Cc\to A$, which is an
$A$-bimodule map, then the following identities hold, for all
$c, d\in \Cc$. 
\begin{equation}\eqlabel{gamma}
c^{(1)}\gamma(c^{(2)}\ot_Ad)= \mu(c\ot_Ad)=
\gamma(c\ot_Ad^{(1)})d^{(2)};\quad
\gamma(c^{(1)}\ot_Ac^{(2)})=\varepsilon(c). 
\end{equation}
The following theorem extends \cite[Theorem 2.6 and Proposition 2.7]{BKW}.

\begin{proposition}\prlabel{coseparable}
Let $\Cc$ be a coseparable coring over a unital ring $A$. Then $\Cc$ is a firm
ring. The categories $\Mm^\Cc$ and $\Mm_\Cc$ are isomorphic, as are the
categories ${}^\Cc\Mm$ and ${}_\Cc\Mm$. Moreover, for all
$P\in\Mm^\Cc$ and $Q\in{^\Cc\Mm}$, the natural morphism 
$$
\xymatrix{
P\ot^\Cc Q \ar[r]^-{\iota}& 
P\ot_A Q \ar[r]^-{\pi}&
P\ot_\Cc Q
} ,
$$
obtained by composing the canonical monomorphism $\iota$ with the canonical
epimorphism $\pi$, is an isomorphism with inverse 
$$
\beta:P\ot_\Cc Q\to P\ot^\Cc Q,\qquad 
\beta(p\ot_\Cc q)= p^{[0]}\ot_A p^{[1]}\cdot q=p \cdot q^{[-1]}\ot_A
q^{[0]}.
$$  
\end{proposition}

\begin{proof}
Take any $M\in\Mm^\Cc$ and define $\mu_M:M\ot_A\Cc \to M$ by the
following composition. 
\[
\mu_M:
\xymatrix{
M\ot_A\Cc \ar[rr]^-{\rho^M\ot_A\Cc} && M\ot_A\Cc\ot_A\Cc
\ar[rr]^-{M\ot_A\gamma} && M\ot_AA\cong M  
} .
\]
Remark that by \equref{gamma}, $\mu_\Cc=\mu$. Let
us check that $\mu_M$ is an associative action, i.e.\ $\mu_M\circ
(\mu_M\ot_A\Cc)=\mu_M\circ (M\ot_A\mu)$, thus in 
particular $\mu$ is an associative multiplication for
$\Cc$. Indeed, for $c,d\in \Cc$,  
\begin{eqnarray*}
\mu_M(m^{[0]}\gamma(m^{[1]}\ot_A c)\ot_Ad) &=&
m^{[0]}\gamma(m^{[1]}\gamma(m^{[2]}\ot_A c)\ot_Ad)\\ 
&=& m^{[0]}\gamma(\gamma(m^{[1]}\ot_A c^{(1)})c^{(2)}\ot_Ad)\\
&=& m^{[0]}\gamma(m^{[1]}\ot_A c^{(1)})\gamma(c^{(2)}\ot_Ad)\\
&=& m^{[0]}\gamma(m^{[1]}\ot_A c^{(1)}\gamma(c^{(2)}\ot_Ad))\\
&=& m^{[0]}\gamma(m^{[1]}\ot_A\mu(c\ot_Ad)) .
\end{eqnarray*}
Next, let us prove that $M$ is a firm $\Cc$-module, that is, the induced map
$\bar{\mu}_M:M\ot_\Cc\Cc\to M$ is an isomorphism with inverse
\[
{\bar{\rho}}^M:
\xymatrix{
M\ar[rr]^-{{\rho}^M} && M\ot_A\Cc \ar[rr]^-{\pi} && M\ot_\Cc\Cc
}.
\]
For all $m\in M$, 
\[
\bar{\mu}_M\circ\bar{\rho}^M(m) 
=\bar{\mu}_M(m^{[0]}\ot_\Cc m^{[1]})
=m^{[0]}\gamma(m^{[1]}\ot_Am^{[2]})
=m^{[0]}\varepsilon(m^{[1]})= m.
\]
On the other hand, for all $m\ot_\Cc c\in M\ot_\Cc\Cc$,  
\begin{eqnarray*}
\bar{\rho}^M\circ\bar{\mu}_M(m\ot_\Cc c)&=&
\bar{\rho}^M(m^{[0]}\gamma(m^{[1]}\ot_Ac))
=m^{[0]}\ot_\Cc m^{[1]}\gamma(m^{[2]}\ot_Ac) \\
&=& m^{[0]}\ot_\Cc \mu(m^{[1]}\ot_Ac) 
=\mu_M(m^{[0]}\ot_Am^{[1]})\ot_\Cc c \\
&=& m^{[0]}\gamma(m^{[1]}\ot_Am^{[2]})\ot_\Cc c = m\ot_\Cc c .
\end{eqnarray*}
This defines a functor $\Xi:\Mm^\Cc\to\Mm_\Cc$ acting on the morphisms as the
identity. This justifies to denote from now on $\mu_M(m\ot_Ac)=m\cdot c$. 
Conversely, take $M\in\Mm_\Cc$. Since $\mu$ is right
$\Cc$-colinear, 
$\Delta_\Cc$ is left $\Cc$-linear. Hence we can define $\rho^M:M\to M\ot_A\Cc$
as 
\[
\rho^M:
\xymatrix{
M \ar[rr]^-{\bar{\mu}_M^{-1}} && M\ot_\Cc\Cc \ar[rr]^-{M\ot_\Cc\Delta_\Cc} &&
M\ot_\Cc\Cc\ot_A\Cc \ar[rr]^-{\bar{\mu}_M\ot_A\Cc} && M\ot_A\Cc  
} .
\]
One can easily check that $(M,\rho^M)\in \Mm^\Cc$, thus we obtain a functor
$\Gamma:\Mm_\Cc\to\Mm^\Cc$ acting on the morphisms as the identity. We
leave it to the reader to verify that $\Xi\circ \Gamma$ and $\Gamma\circ \Xi$
are the identity functors on $\Mm_\Cc$ and $\Mm^\Cc$ respectively. Symmetry
arguments prove ${}^\Cc \Mm\cong {}_\Cc\Mm$.  

To prove the final statement $P\ot_{\Cc} Q\cong P\ot^\Cc Q$, consider the
following diagram, where the row 
and column represent respectively an equalizer and a coequalizer. 
\[
\xymatrix{
&& P\ot_A\Cc\ot_AQ \ar@<-.5ex>[d] \ar@<.5ex>[d] \\
P\ot^\Cc Q \ar[rr]^-\iota && P\ot_A Q \ar@<.5ex>[rr] \ar@<-.5ex>[rr]
\ar[d]^-\pi && P\ot_A\Cc\ot_A Q\\ 
&& P\ot_\Cc Q
}
\]
The map $\alpha:P\ot_AQ\to P\ot_AQ,\
\alpha(p\ot_Aq):=p^{[0]}\ot_Ap^{[1]}\cdot q$ satisfies, for all $p\ot_Aq\in
P\ot_AQ$, 
\begin{eqnarray*}
\alpha(p\ot_A q)&=&
p^{[0]}\ot_Ap^{[1]}\cdot q = p^{[0]}\ot_A
\gamma(p^{[1]}\ot_Aq^{[-1]})q^{[0]}\\ 
&=& p^{[0]}\gamma(p^{[1]}\ot_Aq^{[-1]})\ot_A q^{[0]}
=p\cdot q^{[-1]}\ot_A q^{[0]} .
\end{eqnarray*}
Hence
\begin{eqnarray*}
\big((\rho^P\ot_AQ)\circ \alpha\big)(p\ot_A q)&=&
p^{[0]}\ot_Ap^{[1]}\ot_Ap^{[2]}\cdot q \\
&=&p^{[0]}\ot_Ap^{[1]}\ot_A\gamma(p^{[2]}\ot_Aq^{[-1]}) q^{[0]}\\ 
&=& p^{[0]}\ot_Ap^{[1]}\gamma(p^{[2]}\ot_Aq^{[-1]})\ot_A q^{[0]}\\
&=& p^{[0]}\ot_A\gamma(p^{[1]}\ot_Aq^{[-2]})q^{[-1]}\ot_A q^{[0]}\\
&=& p^{[0]}\gamma(p^{[1]}\ot_Aq^{[-2]})\ot_Aq^{[-1]}\ot_A q^{[0]}\\
&=& p\cdot q^{[-2]}\ot_Aq^{[-1]}\ot_A q^{[0]}
=\big((P\ot_A\rho^Q)\circ\alpha\big)(p\ot_A q),
\end{eqnarray*}
where we used \equref{gamma} in the fourth equation.
Therefore, we obtain by universality of the equalizer a unique morphism
$\alpha':P\ot_AQ\to P\ot^\Cc Q$ such that
$\alpha=\iota\circ\alpha'$. Furthermore, $\alpha'$ is a left
inverse for $\iota$, i.e.\ for all $p\ot_Aq\in P\ot^\Cc Q$, 
\begin{eqnarray*}
(\alpha'\circ\iota)(p\ot_Aq)=p^{[0]}\ot_Ap^{[1]}\cdot q=
p\ot_Aq^{[-1]}\cdot q^{[0]}=p\ot_Aq. 
\end{eqnarray*}
Next we check that
$\alpha'\circ(\mu_P\ot_AQ)=\alpha'\circ(P\ot_A\mu_Q)$. Take any
$p\ot_Ac\ot_Aq\in P\ot_A\Cc\ot_AQ$, then we find 
\begin{eqnarray*}
\alpha(p\cdot c\ot_A q)&=& \alpha(p^{[0]}\gamma(p^{[1]}\ot_Ac)\ot_A q)\\
&=& p^{[0]}\ot_A \gamma(p^{[1]} \gamma(p^{[2]}\ot_Ac)\ot_A q^{[-1]}) q^{[0]}\\
&=& p^{[0]}\ot_A \gamma(\gamma(p^{[1]}\ot_Ac^{(1)})c^{(2)}\ot_A q^{[-1]})
q^{[0]}\\ 
&=& p^{[0]}\ot_A \gamma(p^{[1]}\ot_Ac^{(1)})\gamma(c^{(2)}\ot_A q^{[-1]})
q^{[0]}\\ 
&=& p^{[0]}\ot_A \gamma(p^{[1]}\ot_Ac^{(1)}\gamma(c^{(2)}\ot_A q^{[-1]}))
q^{[0]}\\ 
&=& p^{[0]}\ot_A \gamma(p^{[1]}\ot_A\gamma(c\ot_A q^{[-2]})q^{[-1]})) q^{[0]}\\
&=& p^{[0]}\gamma(p^{[1]}\ot_A\gamma(c\ot_A q^{[-2]})q^{[-1]}))\ot_A  q^{[0]}\\
&=& \alpha(p \ot_A \gamma(c\ot_A q^{[-1]})) q^{[0]})\\
&=& \alpha(p \ot_A c\cdot q) .
\end{eqnarray*}
Since $\alpha$ and $\alpha'$ differ by the monomorphism $\iota$,
from universal property of the coequalizer we therefore obtain a unique
morphism $\beta:P\ot_\Cc Q\to P\ot^\Cc Q$ such that
$\alpha'=\beta\circ\pi$. We then easily compute  
\[
\beta\circ\pi\circ\iota
=\alpha'\circ\iota=P\ot^\Cc Q 
\qquad \textrm{and}\qquad 
\pi\circ\iota\circ\beta\circ\pi
=\pi\circ\iota\circ\alpha'=\pi \circ \alpha =\pi.
\]
Since $\pi$ is an epimorphism, latter identity implies
$\pi\circ\iota\circ\beta=P\ot_\Cc Q$. Hence $\beta$ is an isomorphism with 
inverse $\pi\circ\iota$.   
\end{proof}

As explained in the introduction, the following theorem improves
\cite[Corollary 9.4]{Ver:PhD}, \cite[5.7, 5.8]{Wis:galcom}, \cite[Proposition 
  5.6]{CaDeGrVe:comcor} and is ultimately related to \cite[Theorem
  I]{Schneider90}. Let us emphasize that in the present version of the theorem 
no projectivity condition on the $A$-module $\Cc$ is requested.  

Note that if an $R$-$A$ bimodule $\Sigma$ is an $R$-firmly
  projective right $A$-module then the $R$-bimodule map $R \to \Sigma
  \otimes_A \Sigma^*$, $r\mapsto \sum x_r \otimes_A \xi_r$ in Theorem
  \ref{thm:firm_proj} (iii) induces a (non-unital) ring map $R \to S:= \{\
  x\xi(-)\ |\ x\in \Sigma,\ \xi\in \Sigma^*\ \}\subseteq \End_A(\Sigma)$,
  $r\mapsto \sum x_r \xi_r(-)$.

\begin{theorem}\thlabel{cosepstr}
Let $\Cc$ be a coseparable coring over a unital ring $A$, $R$ a firm ring and
$\Sigma\in{_R\Mm^\Cc}$, such that $\Sigma$ is an $R$-firmly projective right
$A$-module. If $R$ is   
a left ideal in $T=\End^\Cc(\Sigma)$, then the following statements are
equivalent. 

\begin{enumerate}[(i)]
\item $\can:\Sigma^\dagger\ot_R\Sigma\to \Cc$ is surjective;
\item $\can$ is an isomorphism of $A$-corings;
\item $\Hom(\Sigma,-)\ot_RR:\Mm^\Cc\to\Mm_R$ is fully faithful;
\item $-\ot_R\Sigma:\Mm_R\to\Mm^\Cc$ is an equivalence of categories.
\end{enumerate}
\end{theorem}

\begin{proof}
Obviously, $(ii)$ implies $(i)$ and $(vi)$ implies $(iii)$. The implication
$(iii)\Rightarrow(ii)$ 
follows from the structure Theorem for firm Galois comodules, see
\coref{FirmGalois} (vi). We only have to prove that $(i)$ implies $(iv)$. 
Recall from \coref{FirmGalois} (v) that the functor 
$-\ot_R\Sigma:\Mm_R\to\Mm^\Cc$ 
has a right adjoint $-\ot^\Cc\Sigma^\dagger$. Applying
\prref{coseparable}, we obtain the following commutative diagram of functors. 
\begin{equation}\eqlabel{diagramcosep}   
\xymatrix{
\Mm_R \ar@<-0.5ex>[rrrr]_-{-\ot_R\Sigma} \ar@<0.5ex>[rrrrdd]^-{-\ot_R\Sigma}
&&&& \Mm^\Cc \ar@<-0.5ex>[llll]_-{-\ot^\Cc\Sigma^\dagger}
\ar@<0.5ex>[dd]^-{\Xi}\\ 
\\
&&&& \Mm_\Cc \ar@<0.5ex>[uu]^-{\Gamma}
\ar@<0.5ex>[lllluu]^-{-\ot_\Cc\Sigma^\dagger} 
}
\end{equation}
We know from \prref{coseparable} that the vertical functors describe an 
isomorphism of categories, hence the horizontal functors establish an 
equivalence if and only if the diagonal functors do so. The 
diagonal functors are obtained by tensor functors between two module 
categories and can thus be obtained from the Morita context 
$\CC=(R,\Cc,\Sigma,\Sigma^\dagger,\wt,\bt)$, where the connecting maps are 
given by the formulae 
\begin{eqnarray*}
&&\wt: \Sigma\ot_\Cc\Sigma^\dagger\cong \Sigma\ot^\Cc\Sigma^\dagger\cong
\End^\Cc(\Sigma)\ot_RR\cong R;\\ 
&&\bt=\can: \Sigma^\dagger\ot_R\Sigma\to \Cc,\qquad \xi\bt
x=\xi(x^{[0]})x^{[1]} . 
\end{eqnarray*}
Since $R$ is a left ideal in $\End^\Cc(\Sigma)$ by assumption, both
connecting maps of $\CC$ are surjective. Hence 
=the diagonal functors in \equref{diagramcosep} establish an equivalence by
\thref{Kato}, what proves the claim.
\end{proof}

The following example provides another proof for \cite[Corollary
  4.2]{JanTho:descIII} and it also illustrates how our theory goes beyond the
standard case. 

\begin{example}
Let $\iota:B\to A$ be a split extension of unital rings,
i.e.\ such that there exists a $B$-linear morphism $E:A\to B$
such that $E\circ\iota= 
B$. Then the Sweedler coring $\Cc=A\ot_BA$ is coseparable, with
$\mu:A\ot_BA\ot_AA\ot_BA\cong A\ot_BA\ot_BA\to A\ot_BA$ given
by $\mu(a\ot_Ba'\ot_Ba'')=aE(a')\ot_Ba''$ (see
\cite{BKW}). Furthermore, the category $\Mm^\Cc$ is known to be isomorphic to
the category $\Desc(A/B)$ of descent data associated to the ring extension
$\iota$. Note that in this case $\End^\Cc(A)\cong \{\ t\in A\
  |\ t\otimes_B 1_A =1_A \otimes_B t\ \}$ is equal to $B$, hence by
\thref{cosepstr} we find that the categories $\Mm_B$ and $\Desc(A/B)$ are
equivalent by the functor $-\ot_BA$. 

More generally, let $\Sigma\in{_B\Mm_A}$ be a finitely generated projective
right $A$-module which is separable in the sense that the evaluation map
\begin{equation}\label{eq:sep_bim}
\Sigma \ot_A {}_B\Hom(\Sigma,B) \to B,\qquad x \ot_A \xi \mapsto \xi(x) 
\end{equation}
is a split epimorphism of $B$-bimodules. 
Then the associated comatrix coring $\Cc=\Sigma^*\ot_B\Sigma$
is again coseparable (see \cite{BGT}). Since in this case $B
  \to \End_A (\Sigma)$, $b\mapsto (x\mapsto bx)$ is a split extension of
  unital rings, we conclude that $B \cong \End^\Cc(\Sigma)$. Thus we find
that the functor $-\ot_B\Sigma$ is an equivalence between $\Mm_B$ and the
category $\Desc(\Sigma)$ of generalized descent data. 

Consider now a unital ring $A$, a firm ring $R$ and an $R$-firmly
projective right $A$-module $\Sigma$. Assume that $\Sigma$ is a separable
$R$-$A$ bimodule, i.e.\ replacing $B$ by $R$ in \eqref{eq:sep_bim}, we obtain a split
epimorphism of $R$-bimodules. Then the corresponding comatrix $A$-coring
${\mathcal C}:=\Sigma^*\ot_R \Sigma$ is coseparable. Indeed, similarly to
\cite[Theorem 3.5]{BGT}, a bicolinear retraction of the coproduct in
${\mathcal C}$ is given by  
$$
(\Sigma^*\ot_R \Sigma)\ot_A (\Sigma^*\ot_R \Sigma) \to \Sigma^*\ot_R \Sigma,
\qquad 
(\varphi \ot_R z) \ot_A (\psi \ot_R y) \mapsto \sum \varphi\ \xi_r\big( z
\psi(x_r)\big) \ot_R  {^ry},
$$
where 
$r\mapsto \sum x_r
\ot_A \xi_r$ is an $R$-bimodule retraction of \eqref{eq:sep_bim}.

Note that $R$ is a left ideal in $\End^\cC(\Sigma)$. Indeed, taking into
account the explicit form $y\mapsto \sum e_r \ot_A f_r \ot_R \, {}^ry$ of the
${\mathcal C}$-coaction  on $\Sigma$, given in terms of the map 
$R\to \Sigma \ot_A \Sigma^*$, $r\mapsto \sum e_r \ot_A f_r$, encoding $R$-firm
projectivity  of the right $A$-module $\Sigma$, it follows that
$\Phi\in \End_A(\Sigma)$ is a right ${\mathcal C}$ comodule map if and only if  
\begin{equation}\label{eq:colin}
\sum e_r \ot_A f_r \ot_R {}^r \Phi(y)=
\sum \Phi(e_r)\ot_A f_r \ot_R {}^r y,
\end{equation}
for all $y\in \Sigma$. Applying the map
$$
\Sigma\ot_A \Sigma^* \ot_R \Sigma \to \Sigma\qquad
z\ot_A \psi \ot_R y \mapsto \sum \xi_r\big( z \psi(x_r)\big)\, {}^r y
$$
to both sides of \eqref{eq:colin}, we conclude that 
$$
\Phi(y)=\sum (\xi_{r'} \circ \Phi) (r x_{r'})\ {}^{r'}({}^r y)=
\sum (\xi_r \circ \Phi)(x_r)\, {}^ry,
$$
where the last equality follows by the $R$-bilinearity condition $\sum r
x_{r'}\ot_A \xi_{r'} = x_r \ot_A \xi_r r'$. 
This shows that $\Phi r = \xi_r(\Phi(x_r))\in R$, hence $R$ is a left ideal in
$\End^\cC(\Sigma)$, as stated. 

The canonical map corresponding to the $R$-${\mathcal C}$ bicomodule $\Sigma$
is the identity map, so we may conclude by \thref{cosepstr} that $-\ot_R
\Sigma: \Mm_R \to \Mm^{\mathcal C}$ is an equivalence functor. 
\end{example}

\begin{example}
Let $H$ be a Hopf algebra over a commutative ring $k$. Then $H \ot_k H$ admits
the structure of a coseparable $H$-coring cf. \cite[8.8]{BBW:contra}. Right
comodules of the $H$-coring $H\ot_k H$ are known as {\em Hopf modules} of $H$. 
Their category is denoted by $\Mm^H_H$.  

Let $\Sigma$ be an $H$-Hopf module and $T$ be the algebra of Hopf module
endomorphisms of $\Sigma$. Assume that there is a firm ring 
$R$ which is a left ideal in $T$. If $\Sigma$ is an $R$-firmly projective
right $A$-module then, by \thref{cosepstr}, the functor
$-\ot_R\Sigma:\Mm_R\to\Mm^H_H$ is an equivalence if and only if the canonical
map $\can:\Sigma^* \ot_R\Sigma\to H \ot_k H$ is surjective. 

Choose in particular $\Sigma=H$ (with $H$-action given by the multiplication,
and $H$-coaction given by comultiplication) and $R=T\cong k$. Then the inverse
of $\can$ is easily constructed in terms of the antipode of $H$ (see
e.g. \cite[15.5]{BrzWis:cor}). Thus we obtain an alternative proof of the claim
in (iv), what is usually referred to as the {\em Fundamental Theorem of Hopf
  modules}. 
\end{example}

\section*{Acknowledgement}

The first author acknowledges a Bolyai J\'anos Research Scholarship and
financial support of the Hungarian Scientific Research Fund OTKA F67910.

The second author thanks the Fund for Scientific Research--Flanders
(Belgium) (F.W.O.--Vlaanderen) for a Postdoctoral Fellowship.


\begin{thebibliography}{99}
\bibitem{Abuh:cleftentw} J.\ Abuhlail, {\it Morita contexts for corings and
  equivalences,} in: `Hopf algebras in noncommutative geometry and physics'.
  S. Caenepeel and F. Van Oystaeyen (eds.), Marcel Dekker 2005, pp.\ 1--29.
\bibitem{Bass} H.\ Bass, {\it Algebraic K-theory}. Benjamin, New York, 1968.
\bibitem{BDR} M. Beattie, S. D\v asc\v alescu and \c S. Raianu,
{\it Galois extensions for co-Frobenius Hopf algebras},
{J. Algebra} {\bf 198} (1997), 164--183.    
\bibitem{BBW:contra} G. B\"ohm, T. Brzezi\'nski and R. Wisbauer, {\em Monads
  and comonads in module categories}, Preprint
\href{http://arxiv.org/abs/math/0804.1460}{arXiv:0804.1460}. 
\bibitem{BohmVer:cleftex} G. B\"ohm and J. Vercruysse, {\em Morita theory for 
  coring extensions and cleft bicomodules}, Adv. Math. {\bf 209} (2007),
  611--648. {\em Corrigendum}, to be published. See also
\href{http://arxiv.org/abs/math/0601464}{\tt arXiv:math/0601464v2}. 
\bibitem{Brz:corext} T.\ Brzezi\'nski, {\it A note on coring
extensions,} Ann. Univ. Ferrara - Sez. VII - Sc. Mat. 51 (2005), 15--27.
A corrected version is available at 
\href{http://arxiv.org/abs/math/0410020}{\tt arXiv:math/0410020v3}. 
\bibitem{BGT} T. Brzezi\'nski and J. G\'omez Torrecillas,
{\it On comatrix corings and bimodules}, {K-Theory}, {\bf 29} (2003), 101--115.
\bibitem{BKW} T.\ Brzezi{\'n}ski, L.\ Kadison and R.\ Wisbauer,
{\it On coseparable and biseparable corings},
in: `Hopf algebras in noncommutative geometry and physics'.
S. Caenepeel and F. Van Oystaeyen (eds.), Marcel Dekker 2005, pp.\ 71--87.
\bibitem{BrzWis:cor} T.\ Brzezi\'nski and R.\ Wisbauer, {\it Corings and
 Comodules}. Cambridge University Press, Cambridge, 2003.
\bibitem{Cae:bluebook}S.\ Caenepeel, {\it Brauer Groups, Hopf algebras and
  Galois theory}. K-monographs in Mathematics, Kluwer Academic Publishers,
  Dordrecht, 1998. 
\bibitem{CaDeGrVe:comcor} S.\ Caenepeel, E.\ De Groot and J.\
Vercruysse, {\it Galois theory for comatrix corings: descent
theory, Morita theory, Frobenius and separability properties,}
Trans. Amer. Math. Soc. {\bf 359} (2007), 185--226.
\bibitem{CI:semiperfect} S.\ Caenepeel, M.\ Iovanov,
{\it Comodules over semiperfect corings}, in ``Proceedings of the
International Conference on Mathematics and its Applications, ICMA 2004'',
S.L. Kalla and M.M. Chawla (Eds.), Kuwait University, Kuwait, 2005, 135--160. 
\bibitem{CVW:rational} S.\ Caenepeel, J.\ Vercruysse and S.\ Wang, 
{\it Rationality properties for Morita contexts associated to corings}, in
``Hopf algebras in non-commutative geometry and physics", Caenepeel S. and Van
Oystaeyen, F. (eds.), Lect. Notes Pure Appl. Math., Dekker, New York (2005)
113--136. 
\bibitem{CasGT} F.\ Casta\~no-Iglesias and J. G\'omez-Torrecillas, {\it Wide
  Morita contexts}, Comm. Algebra  {\bf 23} (2)  (1995), 601--622.  
\bibitem{ChiDaNa} N. Chifan, S. D\u asc\u alescu and C. N\u ast\u asescu, {\it
  Wide Morita contexts, relative injectivity and equivalence results},
  J. Algebra {\bf 284} (2005), 705--736.
\bibitem{CohFishMont:Morita} M. Cohen, D. Fischman and S.
Montgomery, {\em Hopf Galois extensions, smash products, and
Morita equivalence,} J. Algebra {\bf 133} (1990), 351--372.
\bibitem{Doi:smmor} Y.\ Doi, {\it Generalised smash products and Morita
  contexts for arbitrary Hopf algebras,} in: `Advances in Hopf Algebras'
  J. Bergen and S. Montgomery (eds.), Marcel Dekker 1994.
\bibitem{EK-GT:MoritaDuality} L.\ El Kaoutit and J.\ G\'omez-Torrecillas,
{\it Morita Duality for Corings over Quasi-Frobenius Rings},
in: `Hopf algebras in noncommutative geometry and physics'.
S. Caenepeel and F. Van Oystaeyen (eds.), Marcel Dekker 2005, pp.\ 137--153.
\bibitem{GomVer:ComatFirm} J.\ G\'omez-Torrecillas and J. Vercruysse, {\em
  Comatrix corings and Galois comodules over firm rings},
Algebr. Represent. Theory, {\bf 10} (3) (2007), 271--306.
\bibitem{GraVit:Morita} F.\ Grandjean and E.\ M.\ Vitale, {\it Morita
  equivalence for regular algebras}, Cahiers Topologie G\'eom. Diff\'erentielle
  Cat\'eg. {\bf 39}  (1998), 137--153.  
\bibitem{IV:cofrob} M. Iovanov and J. Vercruysse, {\it Co-Frobenius Corings and
  Related Functors}, J. Pure Appl. Algebra (2008),
  doi:10.1016/j.jpaa.2007.11.015, in press.   
\bibitem{JanTho:descIII} G.\ Janelidze and W.\ Tholen, {\it Facets of descent
  III. Monadic descent for rings and algebras},  Appl. Categ. Structures {\bf
  12}  (2004), 461--477. 
\bibitem{Kato} T.\ Kato and K.\ Ohtake, {\it Morita contexts and equivalences},
  J. Algebra {\bf 61} (1979), 360--366. 
\bibitem{Mar:Morita} L.\ Mar\'in, {\it Morita equivalence based on contexts
  for various categories of modules over associative rings}, J. Pure
  Appl. Algebra  {\bf 133}  (1998), 219--232.  
\bibitem{Mul} B.\ J.\ M\"uller, {\em The quotient category of a Morita
 context}, J. Algebra {\bf 28}  (1974), 389--407. 
\bibitem{Qui:firm} D.\ Quillen, {\it Module theory over nonunital rings},
  unpublished notes, 1997. 
\bibitem{Schneider90}
H.-J.\ Schneider, {\it Principal homogeneous spaces for arbitrary {H}opf
  algebras}, {Israel J. Math.}
{\bf 72}
(1990),
 {167--195}.
\bibitem{Ver:equi} J. Vercruysse, {\em Equivalences between categories of
  modules and categories of comodules}, Acta Math. Sin. (Engl. Ser.) (2008),
  in press. 
\bibitem{Ver:PhD} J. Vercruysse, {\em Galois Theory for Corings and
    Comodules}, PhD thesis. Vrije Universiteit Brussel, 2007.
\bibitem{Wis:com_cat} R.\ Wisbauer, {\em On the category of comodules over
    corings}, in: Mathematics \& Mathematics Education (Bethlehem, 2000),
    pp 325--336, World Sci. Publ., River Edge, NJ, 2002.
\bibitem{Wis:galcom} R.\ Wisbauer, {\em On Galois
comodules}, Comm. Algebra {\bf 34} (2006), 2683-2711.
\end{thebibliography}
\end{document}